\documentclass[secthm,seceqn,amsthm,ussrhead,10pt]{amsart}
\usepackage{amsmath,latexsym}
\usepackage[english]{babel}
\usepackage[psamsfonts]{amssymb}
\usepackage{times}
\usepackage{cite}
\usepackage{pdflscape} 
\usepackage{ulem}
\usepackage[mathcal]{euscript}
\usepackage{tikz}
\usepackage{hyperref}
\usepackage{cancel}
\usetikzlibrary{arrows}

\newtheorem{Th}{Theorem}[section]
\newtheorem{Def}[Th]{Definition}
\newtheorem{Lem}[Th]{Lemma}
\newtheorem{proposition}[Th]{Proposition}
\newtheorem{Remark}[Th]{Remark}
\newtheorem{Example}[Th]{Example}
\newtheorem{corollary}[Th]{Corollary}

\newenvironment{Proof}[1][Proof.]{\begin{trivlist}
\item[\hskip \labelsep {\bfseries #1}]}{\flushright
$\Box$\end{trivlist}}

\begin{document}
	\sloppy

{\Large Complete classification of algebras of level two
\footnote{ The work was supported by RFBR 17-51-04004
the Presidents Programme Support of Young Russian Scientists (MK-2262.2019.1).}
}

\medskip

\medskip

\medskip

\medskip
\textbf{Ivan Kaygorodov$^{a}$, Yury Volkov$^{b}$}
\medskip

{\tiny
$^{a}$ Universidade Federal do ABC, CMCC, Santo Andr\'{e}, Brazil.

$^{b}$ Saint Petersburg state university, Saint Petersburg, Russia.
\smallskip

    E-mail addresses:\smallskip

    Ivan Kaygorodov (kaygorodov.ivan@gmail.com),
    
    Yury Volkov (wolf86\_666@list.ru).

}

       \vspace{0.3cm}

{\bf Abstract.}
The main result of the paper is the classification of all (nonassociative) algebras of level two, i.e. such algebras that maximal chains of nontrivial degenerations starting at them have length two.
During this classification we obtain an estimation of the level of an algebra via its generation type, i.e. the maximal dimension of its one generated subalgebra.
Also we describe all degenerations and levels of algebras of the generation type $1$ with a square zero ideal of codimension $1$.

\smallskip


       \vspace{0.3cm}


       \vspace{0.3cm}

\section{Introduction}

The paper is devote to the classification of algebras of a given level. Algebras in this paper are not assumed to be associative. All algebra structures on a given linear space form an algebraic variety with a natural action of the general linear group.
Orbits under this action correspond to isomorphism classes of algebras. Algebras satisfying some set of polynomial identities constitute a closed subvariety closed under the mentioned action.
There are many papers considering the structure of such subvarieties. One of the main problems in this direction is the description of irreducible components. This problem is called the {\it geometric classification} of algebras.
Examples of a geometric classification in some classes of algebras one can find, for example, in \cite{ikv17,kpv,BC99,S90,GRH,BB09,BB14,kppv}.

Another important notion that is used in the description of varieties of algebras is the {\it degeneration}. One algebra degenerates to another if the closure of the orbit of the first algebra contains the second one.
The description of degenerations helps to describe the irreducible components. For example, if the variety has only finite number of orbits, then any irreducible component is a closure of an orbit of a {\it rigid} algebra, and an algebra is rigid in this case iff there is no nontrivial degeneration to it. On the other hand, degenerations are interesting themselves. There are some papers, where the degeneration graph is constructed for some variety (see, for example, \cite{BB09,gkp17,BB14,BC99,kppv,S90,GRH,kpv,kpv17}).
The notion of a degeneration is closely related to the notions of a contraction and of a deformation.

The notion of the level of an algebra was introduced in \cite{gorb91}. The algebra under consideration has the level $n$ if there is a chain of $n$ nontrivial degenerations that starts at the given algebra and there is no such a chain of length $n+1$. Roughly speaking, the level estimates the complexity of the multiplication of  the given algebra. For example, the unique algebra of the level zero is the algebra with zero multiplication and an algebra has the level one if the closure of its orbit is formed by the zero algebra and the orbit itself. At this moment there are no many results about the levels of algebras. Anticommutative algebras of the first level were classified in \cite{gorb91}, but the classification of all algebras of the first level presented there turned up incorrect. Later the algebras of the first level were classified in \cite{khud13} (see also \cite{IvaPal}). In \cite{gorb93} the author introduced the notion of the infinity level. The infinity level can be expressed in terms of the usual level, and hence the classification of algebras with a given infinity level is much easier than the classification of algebras with a given level. Anticommutative algebras of the second infinity level were classified in \cite{gorb93}. The author made an attempt to classify the anticommutative algebras of the third infinity level in the same paper, but the obtained classification is wrong and can not be taken in account. 
Finally, associative, Lie, Jordan, Leibniz and nilpotent algebras of the level two were classified in \cite{khud15,khud17}.

In the current paper we try to develop a way to classify algebras of small levels. Inspired by the paper \cite{khud13}, we firstly estimate the level of an algebra via its generation type, i.e. the maximal dimension of its one generated subalgebra. We prove that the level of an algebra is not less than its generation type in the case where the generation type is greater of equal to $3$. This estimation is very rough, but is enough for the classification of algebras of small levels. Further we consider different classes of algebras of the generation types one and two and estimate their levels with the help of {\it standard In\"on\"u-Wigner contractions}.

The first type of algebras that we consider is formed by algebras of the generation type $1$ with a square zero ideal of codimension $1$. The anticommutative portion of such algebras was considered in \cite{gorb98}, where they were called almost abelian Lie algebras. Some examples of degenerations between such algebras were given there. In the current paper we describe all degenerations between algebras of the generation type $1$ with a square zero ideal of codimension $1$ and give an explicit formula for the level of such an algebra. Algebras of this type of the first five levels are given in Tables 1-3.

Then we consider algebras of the generation type $1$, whose standard In\"on\"u-Wigner contractions with respect to one dimensional subalgebras have levels not greater than one. It turns out that it is not difficult to classify such algebras. Except the Heisenberg Lie algebras and one algebra of the level $2$, all such algebras have a level not greater than $1$ and moreover have a square zero ideal of codimension $1$. This allows to classify the algebras of the generation type $1$ having the second level.
Note that all anticommutative algebras have the generation type $1$. Thus, we recover the valid part of the results of \cite{gorb93}.

In the remaining part of the paper we consider algebras of the generation type $2$. We give a criterion for a trivial extension of a $2$-dimensional algebra of the generation type $2$ to have the generation type $2$. Then we classify such trivial extensions of the level $2$.
Finally, we consider algebras with an ideal isomorphic to the unique algebra of the generation type $2$ of the first level. All such algebras have degenerations to algebras of a nonantisymmetric bilinear forms. We estimate the level of an algebras of a bilinear form, classify all algebras that have all nonantisymmetric bilinear form degenerations of the level $1$ and estimate the levels of these algebras. In result, we get the classification of all algebras of the level $2$. In particular, we recover and correct the results of \cite{khud15}. It is interesting that all the anticommutative algebras of the second level are Lie algebras and all the alternative algebras of the second level are associative.

\section{Preliminaries}

In this section we introduce some notation and recall some well known definitions and results that we will need in this work.  Note that all the algebras used in this paper are defined in Section \ref{notat} in the end of the paper. In all multiplication tables given there we intend that all omitted products of basic elements are zero. We will be free to use the notation of Section \ref{notat} throughout the paper.

\subsection{Degenerations}
All vector spaces in this paper are over some fixed algebraically closed field ${\bf k}$ and we write simply $dim$, $Hom$ and $\otimes$ instead of $dim_{{\bf k}}$, $Hom_{{\bf k}}$ and $\otimes_{{\bf k}}$.
An algebra in this paper is simply a vector space with a bilinear binary operation. This operation does not have to be associative unlike to the usual notion of an algebra.

Let $V$ be an $n$-dimensional space. Then the {\it set of $n$-dimensional algebra structures} on $V$ is $\mathcal{A}_n=Hom(V\otimes V, V)\cong V^*\otimes V^*\otimes V$. Any $n$-dimensional algebra can be represented by some element of $\mathcal{A}_n$. Two algebras are isomorphic iff they can be represented by the same structure. Moreover, sometimes we will identify a structure from $\mathcal{A}_n$ and an algebra represented by it.
The set $\mathcal{A}_n$ has a structure of the affine variety ${\bf k}^{n^3}$. There is a natural action of the group $GL(V)$ on $\mathcal{A}_n$ defined by the equality $ (g * \mu )(x\otimes y) = g\mu(g^{-1}x\otimes g^{-1}y)$ for $x,y\in V$, $\mu\in \mathcal{A}_2$ and $g\in GL(V)$.
Two structures represent the same algebra iff they belong to the same orbit.

Let $A$ and $B$ be $n$-dimensional algebras. Suppose that $\mu,\chi\in \mathcal{A}_n$ represent $A$ and $B$ respectively. We say that $A$ {\it degenerates} to $B$ and write $A\to B$ if $\chi$ belongs to $\overline{O(\mu)}$. Here, as usually, $O(X)$ denotes the orbit of $X$ and $\overline{X}$ denotes the closer of $X$. We also write $A\not\to B$ if $\chi\not\in\overline{O(\mu)}$. We say that the degeneration $A\to B$ is {\it proper} if $A\not\cong B$. We will write $A\xrightarrow{\not\cong} B$ to emphasize that the degeneration $A\to B$ is proper.
The degeneration $A\xrightarrow{\not\cong} B$ is called {\it primary} if there is no algebra $C$ such that $A\xrightarrow{\not\cong}C$ and $C\xrightarrow{\not\cong} B$.
For an algebra $A$ we introduce the set $A^-=\{B\mid \mbox{there exists a primary degeneration $A\xrightarrow{\not\cong}B$}\}$.

Whenever an $n$-dimensional space named $V$ appears in this paper, we assume that there is some fixed basis $e=(e_1,\dots, e_n)$ of $V$. In this case, for $\mu\in\mathcal{A}_n$, we denote by $\mu_{i,j}^k$ ($1\le i,j,k\le n$) the structure constants of $\mu$ in the basis $e$, i.e. scalars from $\bf k$ such that $\mu(e_i,e_j)=\sum\limits_{k=1}^n\mu_{i,j}^ke_k$. To prove degenerations and nondegenerations we will use the same technique that has been already used in \cite{S90} and \cite{kpv, kppv, kpv17}. In particular, we will be free to use \cite[Lemma 1]{kpv} and facts that easily follow from it.
This lemma asserts the following fact. If $A\to B$, $\mu\in\mathcal{A}_n$ and there is a closed subset $\mathcal{R}\subset\mathcal{A}_n$ invariant under lower triangular transformations of the basis $e_1,\dots,e_n$ such that $\mu\in\mathcal{A}$, then there is a structure $\chi\in\mathcal{R}$ representing $B$. Invariance under lower triangular transformations of the basis $e_1,\dots,e_n$ means that if $\omega\in\mathcal{R}$ and $g\in GL(V)$ has a lower triangular matrix in the basis $e_1,\dots,e_n$, then $g*\omega\in\mathcal{R}$ (see \cite{kpv} for a more detailed discussion).

To prove degenerations, we will use the technique of contractions. Namely, let $\mu,\chi\in \mathcal{A}_n$ represent $A$ and $B$ respectively. Suppose that there are some elements $E_i^t\in V$ ($1\le i\le n$, $t\in{\bf k}^*$) such that $E^t=(E_1^t,\dots,E_n^t)$ is a basis of $V$ for any $t\in{\bf k}^*$ and the structure constants of $\mu$ in this basis are $\mu_{i,j}^k(t)$ for some polynomials $\mu_{i,j}^k(t)\in{\bf k}[t]$. If $\mu_{i,j}^k(0)=\chi_{i,j}^k$ for all $1\le i,j,k\le n$, then $A\to B$. To emphasize that the {\it parametrize basis} $E^t=(E_1^t,\dots,E_n^t)$ ($t\in{\bf k}^*$) gives a degeneration between algebras represented by the structures $\mu$ and $\chi$, we will write $\mu\xrightarrow{E^t}\chi$.
Usually we will simply write down the parametrized basis explicitly above the arrow.

An important role in this paper will be played by a particular case of a degeneration called a {\it standard In\"on\"u-Wigner contraction} (see \cite{IW}). We will call it {\it IW contraction} for short.
Suppose that $A_0$ is an $(n-m)$-dimensional subalgebra of the $n$-dimensional algebra $A$ and $\mu\in\mathcal{A}_n$ is a structure representing $A$ such that $A_0$ corresponds to the subspace $\langle e_{m+1},\dots,e_n\rangle$ of $V$.
Then $\mu\xrightarrow{(te_1,\dots,te_m,e_{m+1},\dots,e_n)}\chi$ for some $\chi\in\mathcal{A}_n$ and the algebra $B$ represented by $\chi$ is called the IW contraction of $A$ with respect to $A_0$. The isomorphism class of the resulting algebra does not depend on the choice of the structure $\mu$ satisfying the condition stated above and always has an ideal $I\subset B$ and a subalgebra $B_0\subset B$ such that $B=B_0\oplus I$ as a vector space, $I^2=0$ and $B_0\cong A_0$ as an algebra. We will call an algebra of such a form a {\it trivial singular extension} of $A_0$ by ${\bf k}^m$.

\subsection{$1$-generated algebras} Let us discuss now some facts about subalgebras generated by one element of an algebra.

\begin{Def}{\rm
Let $A$ be an $n$-dimensional algebra. For $a\in A$, we denote by $A(a)$ the subalgebra of $A$ generated by $a$. The {\it generation type} of $A$ is the dimension of a maximal $1$-generated subalgebra of $A$, i.e. the number $G(A)$ defined by the equality $G(A)=\max\limits_{a\in A}\big(dim\,A(a)\big)$.
}\end{Def}

Let us now choose some structure $\mu\in\mathcal{A}_n$ representing $A$.
It induces a ${\bf k}[x_1,\dots,x_n]$-algebra structure $${\bf k}[x_1,\dots,x_n]\otimes\mu\in Hom\big(({\bf k}[x_1,\dots,x_n]\otimes V)\otimes ({\bf k}[x_1,\dots,x_n]\otimes V), {\bf k}[x_1,\dots,x_n]\otimes V\big)$$ by the equality $({\bf k}[x_1,\dots,x_n]\otimes\mu)\big((f\otimes u)\otimes (g\otimes v)\big)=fg\otimes\mu(u,v)$ for $u,v\in V$, $f,g\in {\bf k}[x_1,\dots,x_n]$.
For two $n$-tuples of polynomials in $n$ variables 
\begin{multline*}
f(x_1,\dots, x_n)=\big(f_1(x_1,\dots, x_n),\dots,f_n(x_1,\dots, x_n)\big)\in ({\bf k}[x_1,\dots,x_n])^n,\\g(x_1,\dots,x_n)=\big(g_1(x_1,\dots,x_n),\dots,g_n(x_1,\dots, x_n)\big)\in ({\bf k}[x_1,\dots,x_n])^n,
\end{multline*}
we define the $n$-tuple $(f\star_{\mu}g)(x_1,\dots,x_n)=\big((f\star_{\mu}g)_1(x_1,\dots,x_n),\dots,(f\star_{\mu}g)_n(x_1,\dots, x_n)\big)\in ({\bf k}[x_1,\dots,x_n])^n$ by the equality
$$({\bf k}[x_1,\dots,x_n]\otimes\mu)\left(\sum\limits_{i=1}^nf_i\otimes e_i,\sum\limits_{i=1}^ng_i\otimes e_i\right)=\sum\limits_{i=1}^n(f\star_{\mu}g)_i\otimes e_i.$$

Let us recall that the Catalan numbers defined by the equality $C_i=\frac{1}{n+1}\binom{2n}{n}$ for $i\ge 0$ satisfy the recursive relation $C_{i}=\sum\limits_{j=0}^{i-1}C_jC_{i-j-1}$. For any $i\ge 1$ we introduce the set of integers $\mathcal{S}_i=\left\{j\left| \sum\limits_{l=0}^{i-2}C_l<j\le \sum\limits_{l=0}^{i-1}C_l \right.\right\}$. It is obvious that $\mathbb{Z}_{>0}$ is a disjoint union of the sets $\mathcal{S}_i$ ($i\ge 1$). Moreover, the formula above guarantees that there exists some bijection $F_{i}:\mathcal{S}_i\rightarrow \bigcup\limits_{j=1}^{i-1}\mathcal{S}_j\times \mathcal{S}_{i-j}$ for any $i\ge 0$. We fix a family of such bijections in this paper and, for $m\in \mathcal{S}_{i}$ we denote by $\mathfrak{l}_m$ and $\mathfrak{r}_m$ such integers that $F_{i}(m)=(\mathfrak{l}_m,\mathfrak{r}_m)$. Thus, for any $m>1$ we have defined two integers $1\le \mathfrak{l}_m,\mathfrak{r}_m<m$. Note that we can choose the bijections $F_i$ ($i\ge 0$) in such a way that $m\le l$ if $\mathfrak{l}_m\le \mathfrak{l}_l$ and $\mathfrak{r}_m\le \mathfrak{r}_l$. We will assume everywhere that the chosen maps $F_i$ satisfy this property.

Now, for a structure $\mu\in\mathcal{A}_n$, we define the $n$-tuples 
$$f^{\mu,i}(x_1,\dots, x_n)=\big(f_1^{\mu,i}(x_1,\dots, x_n),\dots,f_n^{\mu,i}(x_1,\dots, x_n)\big)\in ({\bf k}[x_1,\dots,x_n])^n$$
by induction on $i\ge 1$ in the following way.
Firstly, we define $f^{\mu,1}_j(x_1,\dots,x_n)=x_j$ for all $1\le j\le n$. If $i>1$ and $f^{\mu,j}$ are defined for all $1\le j\le i-1$, then we set $f^{\mu,i}=f^{\mu,\mathfrak{l}_i}\star_{\mu}f^{\mu,\mathfrak{r}_i}$.

Now, it is clear that the vector $v=\sum\limits_{j=1}^n\alpha_je_j\in v$ generates a subalgebra that is generated as a linear space by the vectors $\sum\limits_{j=1}^nf^{i,\mu}_j(\alpha_1,\dots,\alpha_n)e_j$ for $i\ge 1$.

\begin{Def}{\rm
We call $A$ a {\it standard $1$-generated algebra} if there is a direct sum decomposition $A=\oplus_{i\ge 1}A_i$ such that $dim\,A_1=1$, $A_iA_j\subset A_{i+j}$ for all integer numbers $i$ and $j$, and $A$ is generated by $A_1$ as an algebra.
It is easy to see that $G(A)=dim\,A$ for a standard $1$-generated algebra $A$. Note that ${\bf A}_3$ is the unique standard $1$-generated $2$-dimensional algebra structure.
}
\end{Def}

\subsection{Partitions}
In Section \ref{TnGen} we will develop the degenerations of algebras with the generation type $1$ and a square zero ideal of codimension $1$. For this purpose we need the notion of a partition and some facts about it. Note that the notion of a partition was already applied for the studying of the variety of nilpotent matrices in \cite{Ger1}. For more detailed information on partitions we reference the reader to \cite{Ger,BRYL}.

Let us recall that a {\it partition} of the integer number $n$ of the length $l$ is a sequence $a_1,\dots,a_l$ such that $a_1\ge a_2\ge\dots\ge a_l>0$ and $\sum\limits_{i=1}^la_i=n$. In this case we set $len(a)=l$. We denote by $par_n$ the set of all partitions of $n$. We also introduce $par_*=\cup_{n\ge 1}par_n$.
For convenience, we set $a_i=0$ for $i>len(a)$.
Let us define the so-called {\it dominance order} $\succ$ on the set $par_n$. If $a,b\in par_n$, then $a\succeq b$ iff $\sum\limits_{i=1}^ka_i\ge \sum\limits_{i=1}^kb_i$ for all $k\ge 1$. We write $a\succ b$ if $a\succeq b$ and $a\not= b$. It is easy to see that $\succ$ is a partial order on $par_n$. Given $a,b\in par_n$, we say that $b$ is a {\it preceding partition} for $a$ if $a\succ b$ and there is no $c\in par_n$ such that $a\succ c\succ b$. We denote by $a^-$ the set of all preceding partitions for $a$. The following lemma is proved in \cite{BRYL}.

\begin{Lem}\label{succ}
For $a\in par_n$, the set $a^-$ is formed by partitions from the following two sets:
\begin{enumerate}
\item partitions $b$ such that, for some integer $k\ge 1$, $b_k=a_k-1>a_{k+1}$, $b_{k+1}=a_{k+1}+1$, and $b_i=a_i$ if $i\not=k,k+1$;
\item partitions $b$ such that, for some integers $l,k\ge 1$, $b_k=a_k-1=a_{k+1}=\dots=a_{k+l}=a_{k+l+1}+1=b_{k+l+1}$ and $b_i=a_i$ if $i\not=k,k+l+1$.
\end{enumerate}
\end{Lem}

For a partition $a\in par_n$, we will denote by $lev(a)$ the maximal number $m$ such that there exist $a^0,\dots,a^{m-1}\in par_n$ satisfying $a\succ a^{m-1}\succ a^{m-2}\succ\dots\succ a^0$. In other words, $lev(a)$ can be defined by induction in the following way. If $a^-=\varnothing$, then $lev(a)=0$, in the opposite case $lev(a)=1+\max\limits_{b\in a^-}lev(b)$.

Also we will need the sum operation on the set $par_*$. Given partitions $a\in par_n$ and $b\in par_m$, we define their sum $a+b\in par_{n+m}$ by the equality $(a+b)_i=a_i+b_i$ for $i\ge 1$. As usually, the notion of a sum makes sense for any finite family of partitions.

\subsection{Matrices and their full specters.}
Let $\lambda_1,\dots,\lambda_l$ be all the distinct eigenvalues of the matrix $M\in M_n({\bf k})$ and $a^i_1,\dots,a^i_{len(a^i)}$ be the nonincreasing sequence of the sizes of the Jordan blocks corresponding to $\lambda_i$ ($1\le i\le l$) that contains each size as many times as many blocks of the corresponding size $M$ has.
In other words, $a^i\in par_{k_i}$, where $k_i=n-rank(M-\lambda_i\mathcal{E})^n$, and $a^i_1,\dots,a^i_{len(a^i)}$ are such numbers that $rank(M-\lambda_i\mathcal{E})^p=n-\sum_{t=1}^{len(a^i)}\min(a^i_t,p)$ for any $p\ge 1$. Here $\mathcal{E}$ denotes the identity $n\times n$ matrix. Note also that $\sum\limits_{i=1}^lk_i=n$ and $k_i\ge 1$ for any $1\le i\le l$ by definition. We denote the set $\{(\lambda_i,a^i)\}_{i=1}^l\subset {\bf k}\times par_*$ by $FS(M)$ and call it the {\it full specter} of the matrix $M$.

We denote the set $\{FS(M)\mid M\in M_n({\bf k})\}$ of  all possible full specters of $n\times n$ matrices by $FS_n$. The group ${\bf k}^*$ acts on $FS_n$ by the equality $\alpha*\{(\lambda_i,a^i)\}_{i=1}^l=\{(\alpha\lambda_i,a^i)\}_{i=1}^l$ for $\alpha\in{\bf k}^*$ and $\{(\lambda_i,a^i)\}_{i=1}^l\in FS_n$.
It is well known that there is a one to one correspondence between the set $FS_n$ and the set of conjugacy classes of $n\times n$ matrices. It is easy to see also that there is a one to one correspondence between the set $FS_n/{\bf k}^*$ and the set $M_n({\bf k})/\big({\bf k}^*\times GL_n({\bf k})\big)$.
The action of ${\bf k}^*\times GL_n({\bf k})$ on $M_n({\bf k})$ is defined by the equality $(\alpha,U)*M=\alpha UMU^{-1}$ for $\alpha\in{\bf k}^*$, $U\in GL_n({\bf k})$, and $M\in M_n({\bf k})$.
Here and further, for a set $X$ and a group $G$ acting on it, $X/G$ denotes the set of orbits under this action.

Let us introduce, for an integer $l\ge 1$, an $l$-tuple $(i_1,\dots,i_l)$ of nonnegative integers, and $l$-tuple $(\lambda_1,\dots,\lambda_l)$ of elements of ${\bf k}$, the matrix $J_{i_1,\dots,i_l}(\lambda_1,\dots,\lambda_l)$ by the equality
$$
J_{i_1,\dots,i_l}(\lambda_1,\dots,\lambda_l)=
\begin{pmatrix}
J_{i_1}(\lambda_1)&0&\cdots&0&0\\
W_{i_2,i_1}&J_{i_2}(\lambda_2)&&0&0\\
0&W_{i_3,i_2}&\ddots&&\vdots\\
\vdots&0&\ddots&J_{i_{l-1}}(\lambda_{l-1})&0\\
0&\cdots&0&W_{i_l,i_{l-1}}&J_{i_l}(\lambda_l)
\end{pmatrix},$$
$$
\mbox{where }J_i(\lambda)=\left.
\begin{pmatrix}
\lambda&0&0&\cdots&0&0\\
1&\lambda&0&\cdots&0&0\\
0&1&\lambda&\cdots&0&0\\
\vdots&\vdots&\ddots&\ddots&\vdots&\vdots\\
0&0&\cdots&1&\lambda&0\\
0&0&\cdots&0&1&\lambda
\end{pmatrix}\right\}i\mbox{ and }W_{i,j}=\underbrace{\begin{pmatrix}
0&\cdots&0&1\\
0&\cdots&0&0\\
\vdots&0&\vdots&\vdots\\
0&\cdots&0&0
\end{pmatrix}}_{j}\!\!\!\!\!\!\!\!\!\left.\phantom{\begin{pmatrix}\\ \\ \\ \\\end{pmatrix}}\right\}i.
$$
In other words, $J_{i_1,\dots,i_l}(\lambda_1,\dots,\lambda_l)$ is the $\sum\limits_{k=1}^li_k\times \sum\limits_{k=1}^li_k$-matrix that has entries equal to $\lambda_{i_t}$ on the main diagonal from the position $\left(\sum\limits_{k=1}^{t-1}i_k+1\right)$ to the position $\sum\limits_{k=1}^{t}i_k$, has entries equal to $1$ in all the positions on the diagonal that is just below the main diagonal and has all other entries equal to zero. In particular, $J_i(\lambda)$ is the Jordan block of the size $i$ with the eigenvalue $\lambda$.

Let $m$ be some integer and $m_i$ ($1\le i\le k$) be positive integers such that $\sum\limits_{i=1}^km_i=m$. Let $\mathfrak{S}_q$ denote the symmetric group on $q$ elements. Then the group $\mathfrak{S}_{m_1}\times\dots\times \mathfrak{S}_{m_k}$ acts on ${\bf k}^m$ in a natural way. Namely, the $l$-th symmetric group $S_{m_l}$ permutes the components located from the position number $\sum\limits_{i=1}^{l-1}m_i+1$ to the position number $\sum\limits_{i=1}^lm_i$ in the direct product ${\bf k}^m$. Note also that the group ${\bf k^*}$ acts on ${\bf k}^m$ by multiplications. These actions commute and both of them stabilize the zero point. We choose one representative in each orbit of ${\bf k}^m\setminus (0,\dots,0)$ under the action of $\mathfrak{S}_{m_1}\times\dots\times \mathfrak{S}_{m_k}\times {\bf k}^*$ and form a set that we denote by $K^*_{m_1,\dots,m_k}$. Also we choose one representative in each orbit of ${\bf k}^m$ under the action of $\mathfrak{S}_{m_1}\times\dots\times \mathfrak{S}_{m_k}$ and form a set that we denote by $K_{m_1,\dots,m_k}$.

For $b\in par_n$, we introduce $\mathfrak{S}_b=\mathfrak{S}_{b_{len(b)}}\times \mathfrak{S}_{b_{len(b)-1}-b_{len(b)}}\times\dots\times\mathfrak{S}_{b_1-b_2}$. For $b_1$-tuple $(\alpha_1,\dots,\alpha_{b_1}) \in {\bf k}^{b_1}$, let us introduce the matrix
$$
M_{b,\alpha_1,\dots,\alpha_{b_1}}=\begin{pmatrix}
J_{1,\dots,1}(\alpha_1,\dots,\alpha_{b_1})&0&\cdots&0&0\\
0&J_{1,\dots,1}(\alpha_1,\dots,\alpha_{b_2})&\cdots&0&0\\
\vdots&\vdots&\ddots&\vdots&\vdots\\
0&0&\cdots&J_{1,\dots,1}(\alpha_1,\dots,\alpha_{b_{len(b)-1}})&0\\
0&0&\cdots&0&J_{1,\dots,1}(\alpha_1,\dots,\alpha_{b_{len(b)}})
\end{pmatrix}.
$$
If $\lambda_1,\dots,\lambda_l$ are all distinct eigenvalues of the matrix $M_{b,\alpha_1,\dots,\alpha_{b_1}}$, then $FS(M_{b,\alpha_1,\dots,\alpha_{b_1}})=\{(\lambda_i,a^i)\}_{i=1}^l$, where $a^i$ is defined by the equality $a^i_j=\big|\{q\mid 1\le q\le b^i_j,\,\alpha_q=\lambda_i\}\big|$. It is clear that $b=\sum\limits_{i=1}^la^i$ in this case.
It is also not difficult to see that $M_{b,\alpha_1,\dots,\alpha_{b_1}}$ and $M_{b,\beta_1,\dots,\beta_{b_1}}$ are conjugated iff $(\beta_1,\dots,\beta_{b_1})=\sigma(\alpha_{1},\dots,\alpha_{b_1})$ for some $\sigma\in \mathfrak{S}_b$, and $M_{b,\alpha_1,\dots,\alpha_{b_1}}$ and $M_{b,\beta_1,\dots,\beta_{b_1}}$ belong to the same orbit under the action of the group ${\bf k}^*\times GL_n({\bf k})$ iff $(\beta_1,\dots,\beta_{b_1})=\alpha\sigma(\alpha_{1},\dots,\alpha_{b_1})$ for some $\alpha\in{\bf k}^*$ and $\sigma\in \mathfrak{S}_b$.

Let $S=\{(\lambda_i,a^i)\}_{i=1}^l$ be some element of $FS_n$. Let us set $b=\sum\limits_{i=1}^la^i$. It is not difficult to see that one can choose in a unique way $\alpha_1,\dots,\alpha_{b_1}\in\{\lambda_i\}_{1\le i\le l}$ such that $\big|\{q\mid 1\le q\le b^i_j,\,\alpha_q=\lambda_i\}\big|=a^i_j$ and $(\alpha_1,\dots,\alpha_{b_1})\in K_{b_{len(b)},b_{len(b)-1}-b_{len(b)},\dots,b_1-b_2}$. We define $M(S)=M_{b,\alpha_1,\dots,\alpha_{b_1}}$ in this case.
If at least one of the scalars $\lambda_i$ is nonzero, then we also can choose in a unique way $\alpha\in{\bf k}^*$ and $\alpha_1,\dots,\alpha_{b_1}\in\{\lambda_i\}_{1\le i\le l}$ such that $\big|\{q\mid 1\le q\le b^i_j,\,\alpha_q=\lambda_i\}\big|=a^i_j$ and $\alpha(\alpha_1,\dots,\alpha_{b_1})\in K^*_{b_{len(b)},b_{len(b)-1}-b_{len(b)},\dots,b_1-b_2}$. Then we define $M(\bar S)=M_{b,\alpha_1,\dots,\alpha_{b_1}}$, where $\bar S$ denotes the class of $S$ in $FS_n/{\bf k}^*$.  Note that $FS\big(M(S)\big)=S$ and the class of $FS\big(M(\bar S)\big)$ in $FS_n/{\bf k}^*$ equals $\bar S$.

\section{Generation type and level}

In this section we will show that the notions of the generation type and of the level are closely related in the sense that the level of an algebra can be estimated using its generation type. Though the observations of this section are not very surprising, they play a crucial role in our approach to the classification of algebras of low levels. In fact, this approach is inspired by \cite{khud13}, where as the first step of the proof the authors consider the algebras of a generation type more than $1$.

\begin{Lem}\label{closeg}
For any integers $n\ge m\ge 1$, the set $\mathcal{G}_m=\{\mu\in\mathcal{A}_n\mid G(\mu)\le m\}$ is a closed subset of $\mathcal{A}_n$.
\end{Lem}
\begin{Proof}
Let us fix some structure $\mu\in\mathcal{A}_n$. It follows directly from our definitions that $G(\mu)\le m$ iff the rank of the matrix
$$
M_l^{\mu}=
\begin{pmatrix}
f^{\mu,1}_1(\alpha_1,\dots,\alpha_n)&\cdots&f^{\mu,1}_n(\alpha_1,\dots,\alpha_n)\\
\vdots&&\vdots\\
f^{\mu,l}_1(\alpha_1,\dots,\alpha_n)&\cdots&f^{\mu,l}_n(\alpha_1,\dots,\alpha_n)
\end{pmatrix}
$$
is less or equal to $m$ for all $\alpha_1,\dots,\alpha_n\in{\bf k}$ and all $l\ge 0$. It is clear that for a fixed number $l$ this condition is equivalent to some system of polynomial equations in $\mu_{i,j}^k$ ($1\le i,j,k\le n$). Really, the condition $rank(M_l^{\mu})\le m$ is equivalent to the fact that all minors of the dimension $m$ are zero. This gives us a system of polynomial equations in $\alpha_i$ and $\mu_{i,j}^k$ ($1\le i,j,k\le n$). But, since the required equalities have to hold for all $\alpha_1,\dots,\alpha_n\in{\bf k}$, we get polynomial equations in $\mu_{i,j}^k$ ($1\le i,j,k\le n$). Thus, the set
$\mathcal{G}_{m,l}=\{\mu\in\mathcal{A}_n\mid rank(M_l^{\mu})\le m\}$ is closed for any $l\ge 1$. Hence, $\mathcal{G}_m=\bigcap\limits_{l\ge 1}\mathcal{G}_{m,l}$ is closed too.
\end{Proof}

\begin{Def}{\rm
The {\it level} of the $n$-dimensional algebra $A$ is the maximal number $m$ such that there exists a sequence of nontrivial degenerations $A\xrightarrow{\not\cong}A_{m-1}\xrightarrow{\not\cong}\dots\xrightarrow{\not\cong}A_1\xrightarrow{\not\cong}A_0$ for some $n$-dimensional algebras $A_i$ ($0\le i\le m-1$). The level of $A$ is denoted by $lev(A)$.
}\end{Def}

Now we want to find the minimal value of the level of an algebra with a given generation type.
The next lemma shows that standard $1$-generated algebras play a significant role in this problem.

\begin{Lem}\label{AnyToSt}
Suppose that $A$ is an $n$-dimensional algebra with $G(A)=m$. Then $A\to B\oplus {\bf k}^{n-m}$ for some $m$-dimensional standard $1$-generated algebra $B$.
\end{Lem}
\begin{Proof} Let us represent $A$ by a structure $\mu$ such that $\langle e_1,\dots,e_m\rangle$ corresponds to a subalgebra generated by $e_1$. It is clear that $f^{\mu,i}_l(1,0,\dots,0)=0$ for any $i\ge 1$ and any $l>m$ in this case. Let us set $v_i=\sum\limits_{l=1}^mf^{\mu,i}_l(1,0,\dots,0)e_l$ for $i\ge 1$. In particular, $v_1=e_1$. We have $\langle v_i\rangle_{i\ge 1}=\langle e_1,\dots,e_m\rangle$. Then we can choose $1=i_1<i_2<\dots<i_m$ such that $v_{i_l}$ ($1\le l\le m$) are linearly independent and, for any $1\le l<m$ and $i_l<i<i_{l+1}$, the vector $v_i$ belongs to $\langle v_{i_1},\dots,v_{i_l}\rangle$. Let us now choose $d_1,\dots,d_m$ such that $i_l\in \mathcal{S}_{d_l}$ for all $1\le l\le m$. Let us consider the parametrized basis defined by the equalities $E_l^t=t^{d_{l}}v_{i_l}$ for $1\le l\le m$ and $E_l^t=t^{d_{m}}e_l$ for $m+1\le l\le n$.
It is clear from our definitions that, for $1\le k,l\le m$, $\mu(v_{i_k},v_{i_l})=v_s$, where $s=F_{d_k+d_l}^{-1}(i_k,i_l)\in\mathcal{S}_{d_k+d_l}$. Hence, $v_s=\sum\limits_{j=1}^m\alpha_s^jv_{i_j}$ for some $\alpha_s^j\in{\bf k}$ ($1\le j\le m$) such that $\alpha_s^j=0$ if $d_j>d_k+d_l$. Thus, $\mu(E_k^t,E_l^t)=\sum\limits_{j=1}^m\alpha_s^jt^{d_k+d_l-d_j}E_j^t$ and $\mu\xrightarrow{E^t} \chi$ for some $\chi\in\mathcal{A}_n$. It is clear that $\chi(e_i,e_j)=\chi(e_j,e_i)=0$ for $1\le i\le n$ and $m+1\le j\le n$, and $\langle e_1,\dots,e_m\rangle$ is a subalgebra of $\chi$. It remains to show that the restriction of $\chi$ to $\langle e_1,\dots,e_m\rangle$ represents an $m$-dimensional standard $1$-generated algebra. Let us define the grading on $U=\langle e_1,\dots,e_m\rangle$ by the equality $U_d=\langle e_k\mid d_k=d\rangle$.
It is clear from the formula above that $\chi(e_k,e_l)=\sum_{1\le j\le m,d_j=d_k+d_l}\alpha_s^je_j$ for $s=F_{d_k+d_l}^{-1}(i_k,i_l)$, and hence $\chi$ respects the grading on $U$. It is clear that $U_1=\langle e_1\rangle$ is $1$-dimensional, and thus it remains to show that $e_1$ generates $U$ with respect to the structure $\chi$. Let us show using induction on $1\le l\le m$ that $\langle e_1,\dots,e_l\rangle$ lies in the subalgebra generated by $e_1$ with respect to $\chi$. Suppose that the assertion is true for some $l<m$. It is clear that $F_{d_{l+1}}(i_{l+1})=(i_p,i_q)$ for some $1\le p,q\le l$ such that $d_{l+1}=d_p+d_q$.
Really, if it is not so, then $v_{i_{l+1}}=\mu(v_i,v_j)$, where either $v_i\in\langle v_k\rangle_{1\le k< i}$ or $v_j\in\langle v_k\rangle_{1\le k< j}$. Then it follows from the properties of the bijections $F_i$ ($i\ge 0$) that  $v_{i_{l+1}}\in \langle \mu(v_k,v_j)\rangle_{1\le k< i}\subset \langle v_k\rangle_{1\le k< i_{l+1}}$ in the first case and $v_{i_{l+1}}\in \langle \mu(v_i,v_k)\rangle_{1\le k< j}\subset \langle v_k\rangle_{1\le k< i_{l+1}}$ in the second case. This contradicts to the choice of the integers $i_r$ ($1\le r\le m$).
Then we have $v_{i_{l+1}}=\mu(v_{i_p},v_{i_q})$ and $E_{l+1}^t=\mu(E_p^t,E_q^t)$, i.e. $e_{l+1}=\chi(e_p,e_q)$ belongs to the subalgebra generated by $e_1$ with respect to $\chi$. Consequently, the lemma is proved.

\end{Proof}

\begin{corollary}
If $dim\,A=n$, $G(A)=m$ and $lev(A)=\min\limits_{dim A'=n,G(A')=m}lev(A')$, then $A\cong B\oplus {\bf k}^{n-m}$ for some $m$-dimensional standard $1$-generated algebra $B$.
\end{corollary}

The next result estimates the minimal possible level of a standard $1$-generated algebra. This estimation is rough, but it is sufficient for the classification of algebras of low levels.

\begin{Lem}\label{St}
If $A$ is a standard $1$-generated algebra of dimension $n\ge 3$, then $lev(A)\ge n$.
\end{Lem}
\begin{Proof}
By our assumption, $A$ has a grading $A=\oplus_{i\ge 1}A_i$ such that $dim\,A_1=1$ and $A$ is generated by $A_1$ as an algebra. Let us choose a homogeneous basis $e_1,\dots, e_n$ of $A$ such that the degree of $e_i$ is less or equal than the degree of $e_j$ if $i<j$. It is easy to see that $A\xrightarrow{\left(e_1,\dots,e_{n-1},\frac{1}{t}e_n\right)} A'\oplus {\bf k}$ for some standard $1$-generated algebra $A'$ of dimension $n-1$. Thus, it is enough to prove the assertion of the lemma for $n=3$. It is easy to show that any standard $1$-generated algebra of dimension $3$ can be represented either by ${\bf G}$ or by ${\bf G}^{\alpha,\beta}$ for some $(\alpha,\beta)\in K_{1,1}^*$. Since ${\bf G}\xrightarrow{(te_1+te_2, t^2e_2+t^2e_3, t^3e_3)} {\bf G}^{1,1}$, it remains to prove that $lev({\bf G}^{\alpha,\beta})\ge 3$ for any $(\alpha,\beta)\in K_{1,1}^*$.

For the last assertion it is enough to note that there is a sequence of degenerations ${\bf G}^{\alpha,\beta}\xrightarrow{\not\cong}{\bf F}^{\alpha,\beta}\xrightarrow{\not\cong}{\bf A}_3\oplus{\bf k}$. Really, ${\bf G}^{\alpha,\beta}\xrightarrow{(te_1,te_2-te_3,t^2e_3)}{\bf F}^{\alpha,\beta}$, and ${\bf F}^{\alpha,\beta}\xrightarrow{(e_1,e_3,te_2)} {\bf A}_3\oplus{\bf k}$.
\end{Proof}

Lemmas \ref{AnyToSt} and \ref{St} show that, for any algebra $A$ with $G(A)\ge 3$, $lev(A)\ge G(A)$. Moreover $lev(A)\ge G(A)+1$ if $A$ cannot be presented in the form $A=B\oplus {\bf k}^{l}$ for some $G(A)$-dimensional standard $1$-generated algebra $B$.

\begin{Remark}
One can show that $lev(A)\ge 5$ for a standard $1$-generated algebra of dimension $4$. Thus, one can show analogously to the proof of Lemma \ref{St} that $lev(A)\ge n+1$ for a standard $1$-generated algebra of dimension $n\ge 4$. It is interesting to obtain some good estimation for the level of a standard $1$-generated algebra of dimension $n$. For example, it is interesting if this estimation is linear or not.
\end{Remark}


\section{The variety $\mathcal{T}_n$}\label{TnGen}

In this section we introduce the variety $\mathcal{T}_n$ and study its algebraic and geometric properties. This variety is formed by algebras with the generation type $1$ and is important in the study of such algebras, because any IW contraction of an algebra with the generation type $1$ with respect to a $1$-generated subalgebra belongs to $\mathcal{T}_n$.

\subsection{Definition and algebraic description of $\mathcal{T}_n$.} An $n$-dimensional algebra with a square zero ideal of codimension $1$ is an algebra that can be presented by a structure $\mu\in\mathcal{A}_n$ such that
\begin{equation}\label{abel}
\mu(e_i,e_j)=0\mbox{ and }\mu(e_1,e_i),\mu(e_i,e_1)\in \langle e_2,\dots,e_n\rangle\mbox{ for $2\le i,j\le n$}.
\end{equation}

\begin{Lem}\label{G1IC1}
Let $A$ be an algebra presented by the structure $\mu$ satisfying conditions \eqref{abel}. Then $G(A)=1$ iff there exists $\alpha\in{\bf k}$ such that
$$
\mu(e_1,e_1)=\alpha e_1\mbox{ and }\mu(e_1,e_i)+\mu(e_i,e_1)=\alpha e_i\mbox{ for any $2\le i\le n$}.
$$
\end{Lem}
\begin{Proof} It is easy to see that $\mu\big(\sum_{i=1}^n\alpha_ie_i,\sum_{i=1}^n\alpha_ie_i\big)=(\alpha\alpha_1)\sum_{i=1}^n\alpha_ie_i$ for any $\alpha_i\in{\bf k}$ ($1\le i\le n$)  if $\alpha\in{\bf k}$ satisfies the conditions listed in the lemma.
Suppose now that $G(A)=1$. It is clear that $\mu(e_1,e_1)=\alpha e_1$ for some $\alpha\in{\bf k}$ in this case. Since $e_1+e_i$ and $\mu(e_1+e_i,e_1+e_i)=\alpha e_1+\mu(e_1,e_i)+\mu(e_i,e_1)$ have to be linearly dependent, we get $\mu(e_1,e_i)+\mu(e_i,e_1)=\alpha e_i$ for any $2\le i\le n$.
\end{Proof}

Let us denote by $\mathcal{T}_n$ the subset of $\mathcal{A}_n$ formed by structures representing algebras of the generating type $1$ with a square zero ideal of codimension $1$.
It is well known and easy to see that the set of structures representing $n$-dimensional algebras with a square zero ideal of codimension $1$ is a closed subset of $\mathcal{A}_n$. Thus, it follows from Lemma \ref{closeg} that $\mathcal{T}_n$ is a closed subset of $\mathcal{A}_n$.

Let $M=(M_{i,j})_{2\le i,j\le n}$ be an $(n-1)\times(n-1)$ matrix. We define $T_{\alpha}^M\in \mathcal{A}_n$ for $\alpha\in{\bf k}$ in the following way:
$$
T_{\alpha}^M(e_1,e_1)=\alpha e_1,\,\,T_{\alpha}^M(e_1,e_i)=\alpha e_i-T_{\alpha}^M(e_i,e_1)=\sum\limits_{j=2}^nM_{j,i}e_j\,\,(2\le i\le n),\,\,T_{\alpha}^M(e_i,e_j)=0\,\,(2\le i,j\le n).
$$

\begin{corollary}\label{matrixrep} $\mathcal{T}_n=\bigcup\limits_{M\in M_{n-1}({\bf k})}O(T_0^M)\cup \bigcup\limits_{M\in M_{n-1}({\bf k})}O(T_1^M)$. Moreover, $T_r^M$ and $T_s^L$ ($r,s\in\{0,1\}$) lie in the same orbit iff one of the following conditions holds:
\begin{enumerate}
\item $r=s=1$ and there exists $U\in GL_{n-1}({\bf k})$ such that $L=U^{-1}MU$;
\item $r=s=0$ and there exist $U\in GL_{n-1}({\bf k})$ and $\alpha\in{\bf k}^*$ such that $L=\alpha U^{-1}MU$.
\end{enumerate}
\end{corollary}
\begin{Proof} It follows from Lemma \ref{G1IC1} that any structure from $\mathcal{T}_n$ lies in $O(T_{\alpha}^M)$ for some $M\in M_{n-1}({\bf k})$ and $\alpha\in{\bf k}$. If $\alpha\not=0$, then it is easy to see that the structure constants of $T_{\alpha}^M$ in the basis $\frac{e_1}{\alpha},e_2,\dots,e_n$ are the same as the structure constants of $T_1^{\frac{M}{\alpha}}$. Thus, the first assertion is proved.

Suppose that $r,s\in\{0,1\}$. Since the algebra $T_r^M$ is solvable iff $r=0$, it is clear that $T_s^L$ can lie in $O(T_r^M)$ only in the case $r=s$.

Suppose that $T_1^L\in O(T_1^M)$ for some $M,L\in M_{n-1}({\bf k})$. Let $g\in GL(V)$ be such that $g*T_1^L=T_1^M$. Then it is easy to see that the matrix of $g$ in the basis $e_1,\dots,e_n$ has the form
$$
\begin{pmatrix}
1&0&\cdots&0\\
\alpha_2&&&\\
\vdots&&U&\\
\alpha_n&&&
\end{pmatrix}
$$
for some $\alpha_2,\dots,\alpha_n\in{\bf k}$ and $U\in GL_{n-1}({\bf k})$. Then it is easy to see that $T_1^L=g^{-1}*T_1^M=T_1^{U^{-1}MU}$, i.e. $L=U^{-1}MU$.

Finally, if $g*T_0^L=T_0^M$ for some $g\in GL(V)$, then it is easy to see that the matrix of $g$ in the basis $e_1,\dots,e_n$ has the form
$$
\begin{pmatrix}
\alpha&0&\cdots&0\\
\alpha_2&&&\\
\vdots&&U&\\
\alpha_n&&&
\end{pmatrix}
$$
for some $\alpha\in{\bf k}^*$, $\alpha_2,\dots,\alpha_n\in{\bf k}$, and $U\in GL_{n-1}({\bf k})$. Then $T_1^L=g^{-1}*T_1^M=T_1^{\alpha U^{-1}MU}$, i.e. $L=\alpha U^{-1}MU$.
\end{Proof}

We will write simply $T_0^S$ and $T_1^S$ instead of $T_0^{M(\bar S)}$ and $T_1^{M(S)}$ respectively for $S\in FS_{n-1}$. The next corollary follows directly from Corollary \ref{matrixrep}.

\begin{corollary}\label{fsrep} $\mathcal{T}_n=\bigcup\limits_{\overline{S}\in FS_{n-1}/{\bf k}^*}O(T_0^S)\cup \bigcup\limits_{S\in FS_{n-1}}O(T_1^S)$ is a presentation of the variety $\mathcal{T}_n$ as a disjoint union of orbits under the action of $GL_n(V)$.
\end{corollary}

Now we collect some facts about the variety $\mathcal{T}_n$. Namely, we describe its intersections with well known varieties of algebras.
\begin{itemize}
\item The orbits of structures of the form $T_0^S$ ($S\in FS_{n-1}$) correspond to solvable Lie algebras. Such an algebra is nilpotent iff $S=\{(0,a)\}$ for some $a\in par_{n-1}$. The orbits of structures of the form $T_1^S$ ($S\in FS_{n-1}$) are nonsolvable  and not anticommutative, and hence not Lie.
\item Suppose that ${\rm char}{\bf k}\not=2$. Then all nontrivial commutative structures in $\mathcal{T}_n$ belong to the orbit of $T_1^{\left\{\left(\frac{1}{2},(1,\dots,1)\right)\right\}}$. All of these algebras are Jordan and have the level $1$.
\item Associative structures in $\mathcal{T}_n$  belong to the orbits of the structures of one of the forms $$T_0^{\big\{\big(0,(2,\dots, 2,1,\dots,1)\big)\big\}},\,\,T_1^{\big\{\big(1,(1,\dots,1)\big)\big\}},\,\,T_1^{\big\{\big(0,(1,\dots,1)\big)\big\}}\mbox{, and }T_1^{\big\{\big(1,(1,\dots,1)\big),\big(0,(1,\dots,1)\big)\big\}}.$$
\end{itemize}

\subsection{Degenerations in $\mathcal{T}_n$}
The aim of this subsection is to describe all degenerations in the variety $\mathcal{T}_n$. In the end of the subsection we also will give some applications of this description.

\begin{Th}\label{degG1} Suppose that $r,s\in\{0,1\}$, $R,S\in FS_{n-1}$, and $R=\{(\lambda_i,a^i)\}_{i=1}^l$, where $a^i\in par_{k_i}$ and $\lambda_i\in{\bf k}$ for $1\le i\le l$. Then $T_r^R\to T_s^S$ iff one of the following conditions holds:
\begin{enumerate}
\item\label{p1} $r=s=1$ and $S=\{(\lambda_i,b^i)\}_{i=1}^l$ for some $b^i\in par_{k_i}$ ($1\le i\le l$) such that $a^i\succeq b^i$;
\item\label{p2} $r=s=0$ and $S=\{(\alpha\lambda_i,b^i)\}_{i=1}^l$ for some $\alpha\in{\bf k^*}$ and $b^i\in par_{k_i}$ ($1\le i\le l$) such that $a^i\succeq b^i$;
\item\label{p3} $s=0$ and $S=\{(0,b)\}$ for some $b\in par_{n-1}$ such that $\sum\limits_{i=1}^la^i\succeq b$.
\end{enumerate}
\end{Th}
\begin{Proof} For an integer $1\le i\le n-1$, we define $1\le \phi(i)\le l$ as a unique integer such that $0< i-\sum\limits_{j=1}^{\phi(i)-1}k_j\le k_{\phi(i)}$. As before, $\mathcal{E}$ denotes the identity $(n-1)\times (n-1)$ matrix. Firstly, let us show that all the degenerations in $\mathcal{T}_n$ are listed in the theorem.

Let us consider the case $r=1$. Let us set
$$\mathcal{R}=\left\{T_{\alpha}^U\left| \begin{array}{c}\alpha\in{\bf k}, U\in M_{n-1}({\bf k}), U_{i,j}=0\mbox{ for }2\le i<j\le n, U_{i,i}=\alpha\lambda_{\phi(i-1)}\mbox{ for $2\le i\le n$},\mbox{ and}\\
rank\prod\limits_{i=1}^l(U-\alpha\lambda_i\mathcal{E})^{p_i}\le n-1-\sum\limits_{i=1}^l\sum\limits_{j=1}^{len(a^i)}\min(a^i_j,p_i)\mbox{ for any integers $p_i\ge 0$ ($1\le i\le l$)}\end{array}\right.\right\}.$$
It is easy to see that $\mathcal{R}$ is a closed subset of $\mathcal{A}_n$ invariant under lower triangular transformations of the basis $e_1,\dots,e_n$.
Since $O(T_1^R)\cap\mathcal{R}$ is obviously nonempty, any algebra in $\overline{O(T_1^R)}$ can be represented by a structure $T_{\alpha}^U\in\mathcal{R}$.
Thus, if $T_1^R\to T_s^S$, then $T_s^S\in O(T_{\alpha}^U)$ for some $T_{\alpha}^U\in\mathcal{R}$. Let us consider two cases.
\begin{enumerate} 
\item[$\alpha\not=0$.] In this case $T_s^S\cong T_1^{\frac{U}{\alpha}}\in\mathcal{R}$, and hence $s=1$ and $S=FS\left(\frac{U}{\alpha}\right)=\{(\lambda_i,b^i)\}_{i=1}^l$ for some $b^i\in par_{k_i}$. It remains to prove that $a^i\succeq b^i$ for each $1\le i\le l$.
 By the definition of $\mathcal{R}$, we have
\begin{equation}\label{ineq1}
\sum\limits_{j=1}^{len(b^i)}\min(b^i_j,p)=n-1-rank(U-\alpha\lambda_i\mathcal{E})^{p}\ge \sum\limits_{j=1}^{len(a^i)}\min(a^i_j,p)
\end{equation}
for any $p\ge 1$. Let us prove that $\sum\limits_{j=1}^ka^i_j\ge \sum\limits_{j=1}^kb^i_j$ using induction on $k$. Substituting $p=a_1^i$ in \eqref{ineq1}, we get $\sum\limits_{j=1}^{len(b^i)}\min(b^i_j,a^i_1)\ge k_i=\sum\limits_{j=1}^{len(b^i)}b^i_j$. In particular, $b^i_1\le a_1^i$.
Suppose now that  $\sum\limits_{j=1}^{k-1}a^i_j\ge \sum\limits_{j=1}^{k-1}b^i_j$. If $a^i_k\ge b^i_k$, then the induction step is trivial. Let now consider the case $a^i_k< b^i_k$.
Substituting $p=a_k^i$ in \eqref{ineq1}, we get $$k_i-\sum\limits_{j=1}^{k}(b^i_j-a^i_k)\ge \sum\limits_{j=1}^{len(b^i)}\min(b^i_j,a^i_k)\ge \sum\limits_{j=1}^{len(a^i)}\min(a^i_j,a^i_k)=k_i-\sum\limits_{j=1}^{k}(a^i_j-a^i_k),$$ i.e. $\sum\limits_{j=1}^ka^i_j\ge \sum\limits_{j=1}^kb^i_j$.

\item[$\alpha=0$.] In this case, by Corollary \ref{matrixrep}, we have $s=0$ and $S=FS(U)=\{(0,b)\}$ for some $b\in par_{n-1}$. It remains to prove that $\sum\limits_{i=1}^la^i\succeq b$. By definition of $\mathcal{R}$, we have
\begin{equation}\label{ineq2}
\sum\limits_{j=1}^{len(b)}\min\left(b_j,\sum\limits_{i=1}^lp_i\right)=n-1-rank\,U^{\sum\limits_{i=1}^lp_i}\ge \sum\limits_{i=1}^l\sum\limits_{j=1}^{len(a^i)}\min(a^i_j,p_i)
\end{equation}
for any $p_i\ge 0$ ($1\le i\le l$). As before, we are going to prove the inequality $\sum\limits_{j=1}^k\sum\limits_{i=1}^la^i_j\ge \sum\limits_{j=1}^kb_j$ using induction on $k$.
Substituting $p_i=a_1^i$ in \eqref{ineq2} for all $1\le i\le l$, we get $\sum\limits_{j=1}^{len(b)}\min\left(b_j,\sum\limits_{i=1}^la^i_1\right)\ge n-1=\sum\limits_{j=1}^{len(b)}b_j$. In particular, $b_1\le \sum\limits_{i=1}^la_1^i$.
Suppose now that  $\sum\limits_{j=1}^{k-1}\sum\limits_{i=1}^la^i_j\ge \sum\limits_{j=1}^{k-1}b_j$. If $\sum\limits_{i=1}^la^i_k\ge b_k$, then the induction step is trivial. Let now consider the case $\sum\limits_{i=1}^la^i_k< b_k$.
Substituting $p_i=a_k^i$ in \eqref{ineq2} for all $1\le i\le l$, we get $$n-1-\sum\limits_{j=1}^{k}\left(b_j-\sum\limits_{i=1}^la^i_k\right)\ge \sum\limits_{j=1}^{len(b)}\min\left(b_j\sum\limits_{i=1}^l,a^i_k\right)\ge \sum\limits_{i=1}^l\sum\limits_{j=1}^{len(a^i)}\min(a^i_j,a^i_k)=n-1-\sum\limits_{i=1}^l\sum\limits_{j=1}^{k}(a^i_j-a^i_k),$$ i.e. $\sum\limits_{j=1}^k\sum\limits_{i=1}^la^i_j\ge \sum\limits_{j=1}^kb^i_j$.
\end{enumerate}

In the case $r=0$ we set
$$\mathcal{R}=\left\{T_{0}^U\left| \begin{array}{c} U\in M_{n-1}({\bf k}), U_{i,j}=0\mbox{ for }2\le i<j\le n,\mbox{ and $\exists\,\alpha\in{\bf k}$ such that } U_{i,i}=\alpha\lambda_{\phi(i)}\mbox{for $2\le i\le n$ and}\\
rank\prod\limits_{i=1}^l(U-\alpha\lambda_i\mathcal{E})^{p_i}\le n-1-\sum\limits_{i=1}^l\sum\limits_{j=1}^{len(a^i)}\min(a^i_j,p_i)\mbox{ for any integers $p_i\ge 0$ ($1\le i\le l$)}\end{array}\right.\right\}.$$
It is easy to see that $\mathcal{R}$ is a closed subset of $\mathcal{A}_n$ invariant under lower triangular transformations of the basis $e_1,\dots,e_n$. The rest of the proof in the case $r=0$ is analogous to the case $r=1$.

Thus, it remains to show that the degenerations listed in the theorem are valid. It suffices to prove only primary degenerations. According to Lemma \ref{succ} and the statement of the theorem, if $T_r^R\to T_s^S$ is a primary degeneration, then one of the following conditions holds:
\begin{enumerate}
\item $s=r\in\{0,1\}$ and there exist $1\le j\le l$ and $1\le p<q$ such that $S=\{(\lambda_i,b^i)\}_{i=1}^l$, where $b^i_k=a^i_k$ ($1\le i\le l$, $k\ge 1$) if either $i\not=j$ or $k\not\in\{p,q\}$, $b^j_p=a^j_p-1$, and $b^j_q=a^j_q+1$;
\item $s=0$ and $S=\left\{\left(0,\sum\limits_{i=1}^la^i\right)\right\}$.
\end{enumerate}
In the first case we have $T_r^R\cong T_r^M$ and $T_s^S\cong T_r^L$, where
$$
M=\begin{pmatrix}
J_{k-1}(\lambda)&0&0&\cdots&0\\
0&J_{m+1}(\lambda)&0&\cdots&0\\
0&0&&&\\
\vdots&\vdots&&U&\\
0&0&&&
\end{pmatrix}\mbox{ and }
L=\begin{pmatrix}
J_{k}(\lambda)&0&0&\cdots&0\\
0&J_{m}(\lambda)&0&\cdots&0\\
0&0&&&\\
\vdots&\vdots&&U&\\
0&0&&&
\end{pmatrix}
$$
for some $\lambda\in{\bf k}$, integers $1\le k\le m$, and $U\in M_{n-1-k-m}({\bf k})$. Direct calculations show that $T_r^M\xrightarrow{E^t} T_r^L$ for $r\in\{0,1\}$, where
$$
E_i^t=\begin{cases}
e_i,&\mbox{if $i=1$ or $2k+1\le i\le n$},\\
te_{k+1},&\mbox{if $i=2$},\\
te_{i-1},&\mbox{if $3\le i\le k+1$},\\
e_{i}-e_{i-k},&\mbox{if $k+2\le i\le 2k$}.\\
\end{cases}
$$

For the second case, let us introduce $p=\max\limits_{1\le i\le l} len(a^i)$. We assume for simplicity that $len(a^1)\ge len(a^2)\ge\dots\ge len(a^l)$. We have $T_r^R\cong T_r^M$, where
$$
M=\begin{pmatrix}
J_{a^1_1,\dots,a^l_1}(\lambda_1,\dots,\lambda_l)&0&\cdots&0&0\\
0&J_{a^1_2,\dots,a^l_2}(\lambda_1,\dots,\lambda_l)&\cdots&0&0\\
\vdots&\vdots&\ddots&\vdots&\vdots\\
0&0&\cdots&J_{a^1_{p-1},\dots,a^l_{p-1}}(\lambda_1,\dots,\lambda_l)&0\\
0&0&\cdots&0&J_{a^1_p,\dots,a^l_p}(\lambda_1,\dots,\lambda_l)
\end{pmatrix},
$$
and $T_s^S\cong T_0^L$, where
$$
L=\begin{pmatrix}
J_{b_1}(0)&0&\cdots&0&0\\
0&J_{b_2}(0)&\cdots&0&0\\
\vdots&\vdots&\ddots&\vdots&\vdots\\
0&0&\cdots&J_{b_{p-1}}(0)&0\\
0&0&\cdots&0&J_{b_p}(0)
\end{pmatrix},\,\,\,b=\sum\limits_{i=1}^la^i.$$
Direct calculations show that $T_r^M\xrightarrow{E^t} T_0^L$ for $r\in\{0,1\}$, where
$$
E_i^t=\begin{cases}
te_1,&\mbox{if $i=1$},\\
t^{j-1}e_i,&\mbox{if $i=1+\sum\limits_{q=1}^{k-1}b_q+j$, where $1\le k\le p$ and $1\le j\le b_k$}.\\
\end{cases}
$$
\end{Proof}

The next two corollaries follow directly from Theorem \ref{degG1}.

\begin{corollary}\label{primnil} If $a\in par_{n-1}$, then $\left(T_0^{\{(0,a)\}}\right)^-=\left\{T_0^{\{(0,b)\}}\right\}_{b\in a^-}$.
\end{corollary}

\begin{corollary}\label{primgen} Let $r\in\{0,1\}$ and $R=\{(\lambda_i,a^i)\}_{i=1}^l\in FS_{n-1}$. If  either $r=1$ or there exists $1\le i\le l$ with $\lambda_i\not=0$, then
$$
\left(T_r^{R}\right)^-=\left\{T_r^{\{(\lambda_i,a^i)\}_{1\le i\le j-1}\cup\{(\lambda_j,b)\}\cup \{(\lambda_i,a^i)\}_{j+1\le i\le l}}\right\}_{1\le j\le l,b\in (a^j)^-}\cup\left\{T_0^{\left(0,\sum\limits_{i=1}^la^i\right)}\right\}.
$$
\end{corollary}

Now we can compute the level of an algebra from $\mathcal{T}_n$.

\begin{corollary}\label{levelT}
 \begin{enumerate}
\item If $a\in par_{n-1}$, then $lev\left(T_0^{\{(0,a)\}}\right)=lev(a)$.
\item If $r\in\{0,1\}$, $R=\{(\lambda_i,a^i)\}_{i=1}^l\in FS_{n-1}$, and  either $r=1$ or there exists $1\le i\le l$ with $\lambda_i\not=0$, then $lev\left(T_r^{R}\right)=lev\left(\sum\limits_{i=1}^la^i\right)+1$.
\end{enumerate}
\end{corollary}
\begin{Proof} The first assertion follows directly from Corollary \ref{primnil}. Let us prove the second assertion using induction on $\sum\limits_{i=1}^la^i$. If $\left(\sum\limits_{i=1}^la^i\right)^-=\varnothing$, then $(a^i)^-=\varnothing$ and, using Corollary \ref{primgen}, we get
$lev\left(T_r^{R}\right)=1+lev\left(T_0^{\left(0,\sum\limits_{i=1}^la^i\right)}\right)=1$. Suppose that the formula for the level is valid for all $S=\{(\lambda_i,b^i)\}_{i=1}^l\in FS_{n-1}$ such that $\sum\limits_{i=1}^la^i\succ \sum\limits_{i=1}^lb^i$. For integer $1\le j\le l$ and partition $b$ such that $a^j\succ b$, let us introduce $P_{j,b}(a)=\sum\limits_{i=1}^{j-1}a^i+b+\sum\limits_{i=j+1}^{l}a^i\in par_{n-1}$. Note that $\sum\limits_{i=1}^la^i\succ P_{j,b}(a)$. Then, using Corollary \ref{primgen} and the first assertion of this corollary, we get
$$
lev\left(T_r^{R}\right)=1+\max\left(lev\left(\sum\limits_{i=1}^la^i\right),1+\max\limits_{1\le j\le l,b\in (a^j)^-}lev\big(P_{j,b}(a)\big)\right)=lev\left(\sum\limits_{i=1}^la^i\right)+1.
$$
\end{Proof}

\begin{Example}\label{TnNil} Nilpotent algebras of low levels in the variety $\mathcal{T}_n$ are classified in Table 1. All these algebras are of the form $T_0^{\{(0,a)\}}$ for some $a\in par_{n-1}$. This variety is very similar to the variety of $(n-1)\times (n-1)$ nilpotent matrices described in \cite{Ger1}. In particular, Corollary \ref{primnil} and the first part of Corollary \ref{levelT} can be deduced from the just mentioned paper.
\end{Example}

\begin{Example}\label{TnNNil} Suppose that $R=\{(\lambda_i,a^i)\}_{i=1}^l\in FS_{n-1}$ and either $r=1$ or at least one of the scalars $\lambda_i$ is nonzero. Then, by Corollary \ref{levelT}, $lev_n(T_r^R)=m$ iff $lev\left(\sum\limits_{i=1}^la^i\right)=m-1$.
Based on this fact we give classifications of solvable nonnilpotent and nonsolvable algebras of low levels in $\mathcal{T}_n$ in Tables 2 and 3.
\end{Example}

\begin{Remark} It is not difficult to see that $\overline{\cup_{\lambda_1,\dots,\lambda_{n-1}\in{\bf k}}O\left(T_0^{J_{1,\dots,1}(\lambda_1,\dots,\lambda_{n-1})}\right)}$ is an irreducible subvariety of $\mathcal{T}_n$ that coincides with the subvariety of all solvable Lie algebras in $\mathcal{T}_n$.
Moreover, $T_0^{J_{1,\dots,1}(\lambda_1,\dots,\lambda_{n-1})}\in \overline{\cup_{t\in{\bf k}^*}O\left(T_1^{J_{1,\dots,1}\left(\frac{\lambda_1}{t},\dots,\frac{\lambda_{n-1}}{t}\right)}\right)}$, and hence $\mathcal{T}_n=\overline{\cup_{\lambda_1,\dots,\lambda_{n-1}\in{\bf k}}O\left(T_1^{J_{1,\dots,1}(\lambda_1,\dots,\lambda_{n-1})}\right)}$ is an irreducible variety.
\end{Remark}

\section{Classification of algebras of level two}

In this section we classify all the algebras of the level $2$. Note that the described methods can be extended to the study of algebras of higher levels and the obtained results give a reasonable part of the classification of algebras of the level $3$.

\subsection{Algebras of low levels with generation type $1$}
The goal of this subsection is to classify all algebras of the level $2$ with generation type $1$. But firstly let us recall the classification of algebras of the level $1$. For  algebras over the field $\mathbb{C}$ this classification can be found in \cite{khud13}. For algebras over infinite fields the same result can be proved absolutely analogously or can be found in \cite{IvaPal}. However, we give here a short proof for the convenience of the reader.

\begin{proposition}\label{level1} Let $n\ge 2$ be an integer. Then any structure in $\mathcal{A}_n$ corresponding to an algebra of the level $1$ lies in the orbit of exactly one of the structures ${\bf A}_3\oplus{\bf k}^{n-2}$, ${\bf n}_3\oplus{\bf k}^{n-3}$, ${\bf p}^-$ or $\nu^{\alpha}$ ($\alpha\in{\bf k}$).
\end{proposition}
\begin{Proof} It follows from Examples \ref{TnNil} and \ref{TnNNil} that the union of the orbits of ${\bf n}_3\oplus{\bf k}^{n-3}$, ${\bf p}^-$ and $\nu^{\alpha}$ ($\alpha\in{\bf k}$) is exactly the set of all algebras of the level $1$ in $\mathcal{T}_n$. It is easy to see that $\mathcal{R}=\{\mu\in\mathcal{A}_n\mid \mu_{i,j}^k=0\mbox{ for $(i,j,k)\not=(1,1,n)$}\}$
 is a closed subset of $\mathcal{A}_n$ invariant under lower triangular transformations of the basis $e_1,\dots,e_n$. The structure ${\bf A}_3\oplus{\bf k}^{n-2}$ has the level $1$, because its orbit contains all the nontrivial structures from $\mathcal{R}$.

Let now $A$ be an $n$-dimensional algebra. If $G(A)\ge 2$, then it follows from Lemmas \ref{AnyToSt} and \ref{St} that $A\cong B\oplus{\bf k}^{n-2}$ for some $2$-dimensional standard $1$-generated algebra $B$, i.e. $A$ can be represented by ${\bf A}_3\oplus{\bf k}^{n-2}$. If $G(A)=1$ and $A$ has nontrivial multiplication, then there exists $a\in A$ such that $aA+Aa\not=0$. Then the IW contraction with respect to $A(a)$ has to be an algebra of the level $1$ isomorphic to $A$. Since this contraction belongs to $\mathcal{T}_n$, $A$ can be represented by one of the structures ${\bf n}_3\oplus{\bf k}^{n-3}$, ${\bf p}^-$ or $\nu^{\alpha}$ ($\alpha\in{\bf k}$).
\end{Proof}

Let us now prove a lemma about the structure of algebras with the generation type $1$

\begin{Lem}\label{struc}
Let $A$ be an $n$-dimensional algebra with $G(A)=1$. Then $A$ has an $(n-1)$-dimensional subspace $U$ such that $a^2=0$ for all $a\in U$.
\end{Lem}
\begin{Proof} Note that if $a^2=b^2=0$ for $a,b\in A$, then $ab+ba=0$. It is obvious if $a$ and $b$ are linearly dependent. If $a$ and $b$ are linearly independent, then considering $(a+\alpha b)^2=\alpha (ab+ba)$ we obtain that $ab+ba\in \langle a+\alpha b\rangle$ for any $\alpha \in {\bf k}^*$, and hence the equality $ab+ba=0$ holds. If $a,b\in A$ are linearly independent elements such that $a^2\not=0$ and $b^2\not=0$, then rescaling them we may assume that $a^2=a$ and $b^2=b$. Now, considering $(a+\alpha b)^2=a+\alpha (ab+ba)+\alpha^2b=(1+\alpha)(a+\alpha b)+\alpha(ab+ba-a-b)$ we obtain that $ab+ba-a-b\in \langle a+\alpha b\rangle$ for any $\alpha \in {\bf k}^*$, and hence $ab+ba=a+b$. Then $(a-b)^2=0$. Now it is easy to see that we can choose $(n-1)$ linearly independent square zero elements in $A$ that generate a space $U$ with the required properties.
\end{Proof}

We are going to extend the method just used for the classification of algebras of the level $1$ to classify algebras of the level $2$. Thus, we will consider algebras of different generation types separately. In this subsection we consider the algebras of the generation type $1$.
The classification for this case is given in the proposition below.  Note that the $(2m+1)$-dimensional algebra structure $\eta_m$ from Table 6 is known as a {\it nondegenerate Heisenberg Lie algebra}. Note also that $\eta_1={\bf n}_3$ and $\eta_2$ is isomorphic to the $5$-dimensional Lie algebra structure ${\bf g}_1$ that can be found in \cite{GRH,kpv}.

\begin{proposition}\label{GenType1} Let $A$ be an $n$-dimensional algebra with $G(A)=1$.
\begin{enumerate}
\item If $n=2$, then $lev(A)=2$ iff $A$ can be represented by the structure ${\bf E}_4$.
\item If $n=3$, then $lev(A)=2$ iff $A$ can be represented by a structure from the set $$\left\{T_0^{2,\overline{\alpha,\beta}}\right\}_{(\alpha,\beta)\in K_2^*}\cup\left\{T_1^{2,\overline{\alpha,\beta}}\right\}_{(\alpha,\beta)\in K_2}\cup\{{\bf k}\rtimes {\bf E}_4\}.$$
\item If $n=4$, then $lev(A)=2$ iff $A$ can be represented by a structure from the set $$\left\{T_0^{2,\overline{\alpha,\beta}}\right\}_{(\alpha,\beta)\in K_{1,1}^*}\cup\left\{T_1^{2,\overline{\alpha,\beta}}\right\}_{(\alpha,\beta)\in K_{1,1}}\cup\{T_{0}^3,{\bf k}^2\rtimes {\bf E}_4\}.$$
\item If $n\ge 5$, then $lev(A)=2$ iff $A$ can be represented by a structure from the set $$\left\{T_0^{2,\overline{\alpha,\beta}}\right\}_{(\alpha,\beta)\in K_{1,1}^*}\cup\left\{T_1^{2,\overline{\alpha,\beta}}\right\}_{(\alpha,\beta)\in K_{1,1}}\cup\{T_{0}^{2,2},{\bf k}^{n-2}\rtimes {\bf E}_4,\eta_2\oplus{\bf k}^{n-5}\}.$$
\end{enumerate}
\end{proposition}

The rest of this subsection is devoted to the proof of Proposition \ref{GenType1}. Firstly, let us prove that all the algebras mentioned in the proposition have the level $2$. This assertion follows from Examples \ref{TnNil} and \ref{TnNNil} for all the structures, except for the structure ${\bf k}^{n-2}\rtimes {\bf E}_4$ and $\eta_2\oplus{\bf k}^{n-5}$. Thus,  the next lemma finishes the first part of the proof.

\begin{Lem}\label{orbitsoflev2}\begin{enumerate}\item $\overline{O\left({\bf k}^{n-2}\rtimes {\bf E}_4\right)}=O\left(\{{\bf k}^{n-2}\rtimes {\bf E}_4,{\bf p}^-,{\bf k}^n\}\cup\{\nu^{\alpha}\}_{\alpha\in{\bf k}}\right)$ for any $n\ge 2$.
\item $\overline{O\left(\eta_m\oplus{\bf k}^{n-2m-1}\right)}=O\left(\{\eta_l\oplus{\bf k}^{n-2l-1}\}_{1\le l\le m}\cup\{{\bf k}^n\}\right)$ for any $n\ge 2m+1$.
\end{enumerate}
\end{Lem}
\begin{Proof} (1) It is easy to see that ${\bf k}^{n-2}\rtimes {\bf E}_4\xrightarrow{(e_1-e_2,te_2,e_3,\dots,e_n)} {\bf p}^-$, and ${\bf k}^{n-2}\rtimes {\bf E}_4\xrightarrow{\left(\alpha e_1+(1-\alpha)e_2,te_2,e_3,\dots,e_n\right)} \nu^{\alpha}$ for $\alpha\in{\bf k}\setminus\{0\}$.
Let us consider the set
$$
\mathcal{R}=\left\{\mu\in\mathcal{A}_n\left|\begin{array}{c}\mu(e_i,e_j)=0\mbox{ for $2\le i\le n$, $3\le j\le n$},\,\mu(e_i,e_2)=\mu_{2,2}^2e_i,\,\mu(e_i,e_1)=\mu_{2,1}^2e_i\mbox{ for $2\le i\le n$},\\
\mu(e_1,e_i)=\mu_{1,2}^2e_i\mbox{ for $3\le i\le n$},\,\mu(e_1,e_2)=\mu_{2,2}^2e_1+\mu_{1,2}^2e_2,\,\mu(e_1,e_1)=(\mu_{1,2}^2+\mu_{2,1}^2)e_1\end{array}\right.\right\}.
$$
Direct verifications show that $\mathcal{R}$ is a closed subset of $\mathcal{A}_n$ invariant under lower triangular transformations of the basis $e_1,\dots,e_n$. Let us consider $\mu\in \mathcal{R}$. If $\mu_{2,2}^2\not=0$ and $\mu_{1,2}^2\not=0$, then  considering the basis $\frac{e_1}{\mu_{1,2}^2}-\frac{\mu_{2,1}^2e_2}{\mu_{1,2}^2\mu_{2,2}^2},\frac{e_2}{\mu_{2,2}^2},e_3,\dots,e_n$ one can see that $\mu\in O({\bf k}^{n-2}\rtimes {\bf E}_4)$. If  $\mu_{2,2}^2\not=0$ and $\mu_{1,2}^2=0$, then the basis $e_2,e_1-\frac{\mu_{2,1}^2e_2}{\mu_{2,2}^2},e_3,\dots,e_n$ gives $\mu\in O(\nu^0)$. Analogously, it is easy to see that $\mu\in O\left(\nu^{\alpha}\right)$ for $\alpha=\frac{\mu_{1,2}^2}{\mu_{1,2}^2+\mu_{2,1}^2}$ if $\mu_{2,2}^2=0$ and $\mu_{1,2}^2+\mu_{2,1}^2\not=0$, and $\mu\in O\left({\bf p}^-\right)$ if $\mu_{2,2}^2=0$, $\mu_{1,2}^2+\mu_{2,1}^2=0$ and $\mu_{1,2}^2\not=0$.

(2) Note that $\eta_m\xrightarrow{(e_1,\dots,e_{2m-2},e_{2m+1},e_{2m},te_{2m-1})} \eta_{m-1}\oplus{\bf k}^2$. On the other hand, a nilpotent Lie algebra $A$ with $dim\,A^2=1$ can be represented by $\eta_l\oplus{\bf k}^{n-2l-1}$ for some $1\le l<\frac{n}{2}$. Thus, $\overline{O\left(\eta_m\oplus{\bf k}^{n-2m-1}\right)}$ contains only the orbits listed in the statement of the lemma.
\end{Proof}

It remains to prove that any algebra with the level $2$ and the generation type $1$ can be represented by a structure from Proposition \ref{GenType1}.
The main idea of the proof is the following. If the algebra $A$ has generation type $1$, then any IW contraction of $A$ with respect to a $1$-generated subalgebra is an algebra from $\mathcal{T}_n$.
Thus, the main step of our proof is the classification of algebras $A$ with $lev(A)=2$ and $G(A)=1$ such that any IW contraction of $A$ with respect to a $1$-dimensional subalgebra has a level not greater than $1$.
Suppose that $A$ satisfies the just stated conditions.
We say that $a\in A\setminus\{0\}$ is of {\it $X$-type}, where $X\in\{{\bf n}_3\oplus{\bf k}^{n-3},{\bf p}^-\}\cup\{\nu^{\alpha}\}_{\alpha\in{\bf k}}$, if the IW contraction of $A$ with respect to $A(a)$ can be represented by the structure $X$. If the corresponding IW contraction is trivial, we say that $a$ is of {\it $0$-type}. We will also write simply ${\bf n}_3$-type instead of ${\bf n}_3\oplus{\bf k}^{n-3}$-type. By Proposition \ref{level1}, any $a\in A$ is of $0$-type, ${\bf n}_3$-type, ${\bf p}^-$-type or $\nu^{\alpha}$-type for some $\alpha\in{\bf k}$.

\begin{Lem}\label{nu_p_lem}
Let $A$ be an algebra with $G(A)=1$ such that any element of $A$ has $0$-type, ${\bf p}^-$-type or $\nu^{\alpha}$-type.
Suppose that $a,b\in A$ are two linearly independent elements. Then
\begin{enumerate}
\item if $a$ and $b$ are of ${\bf p}^-$-type, then $b-\gamma a$ is of $0$-type for some $\gamma\in{\bf k}$;
\item if $a$ is of $\nu^{\alpha}$-type for some $\alpha\in{\bf k}$ and $b$ is of ${\bf p}^-$-type, then $a-\gamma b$ is of $\nu^1$-type for some $\gamma\in{\bf k}$.
\end{enumerate}
\end{Lem}
\begin{Proof} (1) Suppose that $a$ and $b$ are of ${\bf p}^-$-type. After a rescaling by nonzero scalars, we may assume that $ac=c+f_1(c)a$, $ca=-c+f_2(c)a$, $bc=c+g_1(c)b$, $cb=-c+g_2(c)b$ for any $c\in A$, where $f_1,f_2,g_1,g_2\in Hom(A,{\bf k})$ satisfy the equality $-f_1(a)=f_2(a)=-g_1(b)=g_2(b)=1$.
In particular, $a(b-a)=(a-b)a=b-a$. Since $a$ and $b-a$ are linearly independent, it is easy to see that $b-a$ cannot have ${\bf p}^-$-type or $\nu^{\alpha}$-type for some $\alpha\in{\bf k}$. Thus, $b-a$ is of $0$-type.

(2) After a rescaling by nonzero scalars, we may assume that $ac=\alpha c+f_1(c)a$, $ca=(1-\alpha)c+f_2(c)a$, $bc=c+g_1(c)b$, $cb=-c+g_2(c)b$ for any $c\in A$, where $f_1,f_2,g_1,g_2\in Hom(A,{\bf k})$ satisfy the equalities $f_1(a)=1-\alpha$, $f_2(a)=\alpha$, $-g_1(b)=g_2(b)=1$.
In particular, $\big(a+(1-\alpha)b\big)b=\alpha b-a=-\big(a+(1-\alpha)b\big)+b$, $b\big(a+(1-\alpha)b\big)=a+(1-\alpha)b$ and, hence, $a+(1-\alpha)b$ is of $\nu^1$-type.
\end{Proof}

\begin{corollary}\label{nu_ptype} Suppose that $G(A)=1$ and any IW contraction of the $n$-dimensional algebra $A$ with respect to a $1$-dimensional subalgebra has a level not greater than $1$. If $A$ does not have an element of ${\bf n}_3$-type, then it either has a level not greater than $1$ or can be represented by the structure ${\bf k}^{n-2}\rtimes {\bf E}_4$.
\end{corollary}
\begin{Proof} It easily follows from Lemmas \ref{struc} and \ref{nu_p_lem} that $A$ has a basis $a_1,\dots,a_n$ such that $a_i$ is of $0$-type for $i\ge 3$ and either $a_2$ is of $0$-type too or $a_1$ is of $\nu^1$-type and $a_2$ is of ${\bf p}^-$-type. If $a_2$ is of $0$-type, then it is clear that $A$ is either trivial or of the level $1$. In the second case, after rescaling of the elements $a_1$ and $a_2$, we have $a_1a_1=a_1$, $a_1a_2=-a_1+a_2$, $a_2a_1=a_1$, $a_1a_i=a_2a_i=-a_ia_2=a_i$ for $3\le i\le n$, and all the remaining products of basic elements equal to zero. Changing $a_2$ by $a_1-a_2$, one can see that $A$ can be represented by ${\bf k}^{n-2}\rtimes {\bf E}_4$.
\end{Proof}

\begin{Lem}\label{n3type} Suppose that $G(A)=1$ and any IW contraction of the algebra $A$ with respect to a $1$-dimensional subalgebra has a level not greater than $1$. If $A$ has an element of ${\bf n}_3$-type, then it can be represented by the structure $\eta_m\oplus{\bf k}^{n-2m-1}$ for some $1\le m<\frac{n}{2}$.
\end{Lem}
\begin{Proof} Let us represent $A$ by a structure $\mu\in\mathcal{A}_n$ such that the IW contraction of $\mu$ with respect to the subalgebra generated by $e_1$ equals ${\bf n}_3\oplus {\bf k}^{n-3}$. This means that $\mu(e_1,e_1)=0$, $\mu(e_1,e_2)=\mu_{1,2}^1e_1+e_3$, $\mu(e_2,e_1)=\mu_{2,1}^1e_1-e_3$, $\mu(e_1,e_i)=\mu_{1,i}^1e_1$ and $\mu(e_i,e_1)=\mu_{i,1}^1e_1$ for $3\le i\le n$. Changing $e_3$ by $e_3+\mu_{1,2}e_1$ we may assume that $\mu_{1,2}^1=0$. Since $\mu(e_1,e_2)=e_3$, the basic element $e_2$ cannot be of $0$-type, ${\bf p}^-$-type or $\nu^{\alpha}$-type for any $\alpha\in{\bf k}$. Thus, $e_2$ is of ${\bf n}_3$-type, and hence $e_1e_2+e_2e_1\subset\langle e_2\rangle$, i.e. $\mu_{2,1}^1=\mu_{2,3}^3=\mu_{3,2}^3=0$. Subtracting $\mu_{2,i}^3e_1$ from $e_i$ for $i\ge 4$, we may assume that
\begin{multline*}
\mu(e_1,e_1)=\mu(e_2,e_2)=0,\,\mu(e_1,e_2)=-\mu(e_2,e_1)=e_3,\\\mu(e_1,e_i)=\mu_{1,i}^1e_1,\,\mu(e_i,e_1)=\mu_{i,1}^1e_1,\,\mu(e_2,e_i)=\mu_{2,i}^2e_2,\,\mu(e_i,e_2)=\mu_{i,2}^2e_2\,\,(3\le i\le n).
\end{multline*}

Suppose that $e_3$ is of $p^-$-type. After rescaling $e_1$ and $e_3$, we may assume that $\mu(e_3,e_2)=e_2$. Then we have $\mu(e_1+e_3,e_2)=e_2+e_3$ and $\mu(e_1+e_3,e_3)=-(e_1+e_3)+e_3$. Thus, the IW contraction of $\mu$  with respect to the subalgebra generated by $e_1+e_3$ is not of the first level. The obtained contradiction shows that $e_3$ cannot be of $p^-$-type. Analogously, $e_3$ cannot be of $\nu^{\alpha}$-type for any $\alpha\in{\bf k}$. Then $e_3$ is either of $0$-type or of ${\bf n}_3$-type and, in particular, $\mu(e_1,e_3)=\mu(e_3,e_1)=\mu(e_2,e_3)=\mu(e_3,e_2)=0$.
Now, for $i\ge 4$, the argument as above shows that if $e_i$ is of ${\bf p}^-$-type or of $\nu^{\alpha}$-type for some $\alpha\in{\bf k}$, then the IW contraction of $\mu$  with respect to the subalgebra generated by $e_1+e_i$ is not of the first level. Thus, any element of $A$ is either of $0$-type or of ${\bf n}_3$-type, and thus $\mu_{1,i}^1=\mu_{i,1}^1=\mu_{2,i}^2=\mu_{i,2}^2=0$ for any $i\ge 3$.
In particular, $a^2=0$ for any $a\in A$, i.e. $A$ is anticommutative.

Suppose that $e_3$ is of ${\bf n}_3$-type. Then we may assume that $e_4e_3\not=0$. Since $(e_1+e_4)e_2=e_3$, the element $e_1+e_4$ has to be of ${\bf n}_3$-type, and hence $e_4e_3=(e_1+e_4)e_3=\alpha (e_1+e_4)$ for some $\alpha\in {\bf k}$. Since $e_3$ is of ${\bf n}_3$-type, we have $\alpha=0$. The obtained contradiction shows that $e_3$ has to be of $0$-type. Moreover, since any element of $A$ is either of $0$-type or of ${\bf n}_3$-type, we have $\mu(e_3,e_i)=\mu(e_i,e_3)=0$ for all $1\le i\le n$. Suppose that there exist $4\le i,j\le n$ such that $\mu(e_i,e_j)\not\in \langle e_3\rangle$. We may assume that $i=4$ and $j=5$.
Since $(e_1+e_4)e_2=e_3$ and $e_1+e_4$ has to be of ${\bf n}_3$-type, we have $e_4e_5=(e_1+e_4)e_5=\alpha_1 (e_1+e_4)+\beta_1 e_3$ for some $\alpha_1\in {\bf k}^*$ and $\beta_1\in {\bf k}$. Analogously, considering $e_2+e_5$, we get $e_4e_5=e_4(e_2+e_5)=\alpha_2 (e_2+e_5)+\beta_2 e_3$ for some $\alpha_2\in {\bf k}^*$ and $\beta_2\in {\bf k}$. This contradicts the linearly independence of $e_1+e_4$, $e_2+e_5$ and $e_3$. Thus, $\mu(e_i,e_j)\in \langle e_3\rangle$ for any $1\le i,j\le n$. Thus, $A$ is an anticommutative nilpotent algebra with $dim\,A^2=1$. The statement of the lemma easily follows from this fact.
\end{Proof}

\begin{Proof}[Proof of Proposition \ref{GenType1}.]
It remains to prove that any algebra of the level $2$ can be represented by a structure from the statement of the proposition.
Suppose that the $n$-dimensional algebra $A$ has the level $2$. 
If there exists $a\in A$ such that the IW contraction of $A$ with respect to $A(a)$ has a level greater than $1$, then this contraction has the level $2$ and it is isomorphic to $A$. In particular, $A$ belongs to $\mathcal{T}_n$ in this case and the required assertion follows from Examples \ref{TnNil} and \ref{TnNNil}.
If the level of any IW contraction of $A$ with respect  a $1$-dimensional subalgebra has a level not greater than $1$, then the required assertion follows from Proposition \ref{level1}, Corollary \ref{nu_ptype} and Lemma \ref{n3type}.
\end{Proof}

Note that Corollary \ref{nu_ptype} and Lemmas \ref{orbitsoflev2} and \ref{n3type} give the following interesting result that can be useful for the classification of algebras of levels higher than $2$.

\begin{corollary} If $lev(A)=m>2$, $G(A)=1$ and any IW contraction of the algebra $A$ with respect to a $1$-dimensional subalgebra has a level not greater than $1$, then $n\ge 2m+1$ and $A$ can be represented by the structure $\eta_m\oplus{\bf k}^{n-2m-1}$.
\end{corollary}

Since any anticommutative algebra by definition has the generation type $1$, we get the following classification of anticommutative algebras of the level $2$.

\begin{corollary}\label{levelLie} Let $A$ be an $n$-dimensional anticommutative algebra.
\begin{enumerate}
\item If $n=2$, then $lev(A)\le 1$.
\item If $n=3$, then $lev(A)=2$ iff $A$ can be represented by the structure $T_0^{2,\overline{\alpha,\beta}}$ for some $(\alpha,\beta)\in K_2^*$.
\item If $n=4$, then $lev(A)=2$ iff $A$ can be represented by a structure from the set $\left\{T_0^{2,\overline{\alpha,\beta}}\right\}_{(\alpha,\beta)\in K_{1,1}^*}\cup\{T_{0}^3\}.$
\item If $n\ge 5$, then $lev(A)=2$ iff $A$ can be represented by a structure from the set $\left\{T_0^{2,\overline{\alpha,\beta}}\right\}_{(\alpha,\beta)\in K_{1,1}^*}\cup\{T_{0}^{2,2},\eta_2\oplus{\bf k}^{n-5}\}.$
\end{enumerate}
In particular, all anticommutative algebras of the level $2$ are Lie algebras.
\end{corollary}

\begin{Remark} One can check that the results of Corollary \ref{levelLie} give the same classification as \cite[Proposition 3.3, Theorem 3.5]{khud15} modulo some mistakes in the mentioned paper.
Namely, the author of \cite{khud15} used a wrong version of the description of the degenerations of $3$-dimensional Lie algebras and gave a wrong classification for the case $n=3$.  Also, there is a misprint in \cite{khud15} excluding the algebra isomorphic to $T_0^{2,\overline{0,1}}$ from the classification in the case $n\ge 5$.
\end{Remark}

\subsection{Extensions of $2$-dimensional algebras with generation type $2$}
In the studying of levels of algebras with generation type $2$, we are going to use the same tool as in the case of generation type $1$, i.e. IW contractions. In the case of generation type $2$, we are going to apply them with respect to $2$-dimensional $1$-generated algebras. This subsection is devoted to the algebras that can be obtained in result, i.e. to trivial singular extensions with generation type $2$ of $2$-dimensional algebras with generation type $2$.

Let $C$ be a $2$-dimensional algebra with $G(C)=2$. A trivial singular extension of $C$ is an $n$-dimensional algebra $A$ that has an ideal $I\subset A$ and an injective algebra homomorphism $\phi:C\rightarrow A$ such that $I^2=0$ and $A=\phi(C)\oplus I$ as a vector space.
It follows from the results of \cite{kpv17} that $C$ has an element $a$ such that $a$ and $a^2$ are linearly dependent. Then it is easy to show that $C$ can be represented by a structure $\chi\in\mathcal{A}_2$ such that $\chi_{1,1}^2=1$ and $\chi(e_2,e_2)=\chi_{2,2}^2e_2$, where $\chi_{2,2}^2\in\{0,1\}$. We will denote the set of such structures by $\mathcal{\tilde A}_2$.
Suppose that $A$, $I$ and $\phi:C\rightarrow A$ are as above. Let us represent the algebra $A$ by a structure $\mu\in\mathcal{A}_n$ such that $\langle e_3,\dots,e_n\rangle$ corresponds to the ideal $I$, $\langle e_1,e_2\rangle$ corresponds to the subalgebra $\phi(C)$ and moreover $\mu_{i,j}^k=\chi_{i,j}^k$ for $1\le i,j,k\le 2$. The structure $\mu\in\mathcal{A}_n$ is fully determined by the structure $\chi\in\mathcal{A}_2$ and four matrices $L_1,R_1,L_2,R_2\in M_{n-2}({\bf k})$ such that $\mu_{i,j}^k=(L_i)_{kj}$ and $\mu_{j,i}^k=(R_i)_{kj}$ for $i=1,2$ and $3\le j,k\le n$. Here, for the convenience, we enumerate the rows and the columns of all the $(n-2)\times (n-2)$ matrices under consideration by the numbers from $3$ to $n$. Moreover, we identify these matrices with the corresponding linear transformations of $\langle e_3,\dots,e_n\rangle$. We will denote by ${\bf k}^{n-2}\rtimes_{(L_1,R_1,L_2,R_2)}\chi$ the structure determined by the structure $\chi$ and the matrices $L_1,R_1,L_2,R_2$. Let $\mathcal{E}$ denote the $(n-2)\times (n-2)$ identity matrix and $S$ denote the matrix $L_1+R_1-\chi_{1,1}^1\mathcal{E}$.

\begin{proposition}\label{G2eq}
If $\mu={\bf k}^{n-2}\rtimes_{(L_1,R_1,L_2,R_2)}\chi$ for some $\chi\in\mathcal{\tilde A}_2$ and $L_1,R_1,L_2,R_2\in M_{n-2}({\bf k})$, then $G(\mu)=2$ iff
$$
L_2=(\chi_{2,1}^2\mathcal{E}-R_1)S+\chi_{2,1}^1\mathcal{E},\,R_2=(\chi_{1,2}^2\mathcal{E}-L_1)S+\chi_{1,2}^1\mathcal{E}\mbox{ and }S^3+(\chi_{1,1}^1-\chi_{1,2}^2-\chi_{2,1}^2)S^2+(\chi_{2,2}^2-\chi_{1,2}^1-\chi_{2,1}^1)S=0.
$$
\end{proposition}
\begin{Proof} It is clear that $\mu$ represents an algebra with the generating type $2$ iff, for any $v\in V$, we have $\mu(v,w),\mu(w,v),\mu(w,w)\in\langle v,w\rangle$, where $w=\mu(v,v)$. Since $\mu(v,v)=0$ for any $v\in\langle e_3,\dots,e_n\rangle$, we have to check the required condition for elements of the form $e_1+te_2+v$ and $e_2+v$, where $t\in{\bf k}$ and $v\in\langle e_3,\dots,e_n\rangle$. Let us introduce $\alpha_t=1+(\chi_{1,2}^2+\chi_{2,1}^2-\chi_{1,1}^1)t+(\chi_{2,2}^2-\chi_{1,2}^1-\chi_{2,1}^1)t^2$ and $M_t=S+t\big(L_2+R_2-(\chi_{1,2}^1+\chi_{2,1}^1)\mathcal{E}\big)$. For $u_t=e_1+te_2+v$, direct calculations show that
$w_t=\mu(u_t,u_t)-\big(\chi_{1,1}^1+(\chi_{1,2}^1+\chi_{2,1}^1)t\big)u_t=\alpha_t e_2+M_t v.$ Now we have to check that $\mu(u_t,w_t),\mu(w_t,u_t),\mu(w_t,w_t)\in\langle u_t,w_t\rangle$. We have
$$
\mu(u_t,w_t)-\chi_{1,2}^1\alpha_tu_t-\big(\chi_{1,2}^2+(\chi_{2,2}^2-\chi_{1,2}^1)t\big)w_t=\big((L_1-\chi_{1,2}^2\mathcal{E})M_t+t(L_2-\chi_{2,2}^2\mathcal{E}+\chi_{1,2}^1\mathcal{E})M_t+\alpha_t(R_2-\chi_{1,2}^1\mathcal{E})\big)v.
$$
It is clear that the obtained vector lies in $\langle u_t,w_t\rangle$ for any $t\in{\bf k}$ and any $v\in\langle e_3,\dots,e_n\rangle$ iff $$(L_1-\chi_{1,2}^2\mathcal{E})M_t+t(L_2-\chi_{2,2}^2\mathcal{E}+\chi_{1,2}^1\mathcal{E})M_t+\alpha_t(R_2-\chi_{1,2}^1\mathcal{E})=0$$ for any $t$. Considering the coefficients at the zero, first and second degrees of $t$, we get the equalities
\begin{multline*}
(L_1-\chi_{1,2}^2\mathcal{E})S+R_2-\chi_{1,2}^1\mathcal{E}=0,\\
(L_1-\chi_{1,2}^2\mathcal{E})\big(L_2+R_2-(\chi_{1,2}^1+\chi_{2,1}^1)\mathcal{E}\big)+(L_2-\chi_{2,2}^2\mathcal{E}+\chi_{1,2}^1\mathcal{E})S+(\chi_{1,2}^2+\chi_{2,1}^2-\chi_{1,1}^1)(R_2-\chi_{1,2}^1\mathcal{E})=0,\\
(L_2-\chi_{2,2}^2\mathcal{E}+\chi_{1,2}^1\mathcal{E})\big(L_2+R_2-(\chi_{1,2}^1+\chi_{2,1}^1)\mathcal{E}\big)+(\chi_{2,2}^2-\chi_{1,2}^1-\chi_{2,1}^1)(R_2-\chi_{1,2}^1\mathcal{E})=0.
\end{multline*}
Thus, the formula for $R_2$ stated in the proposition has to be satisfied. Analogously, the formula for $L_2$ follows from $\mu(w_t,u_t)\in\langle u_t,w_t\rangle$. Substituting the obtained values of $L_2$ and $R_2$ in the second of the obtained equalities, we get the last required equality.
On the other hand, if all the required equalities are satisfied, then direct calculations show that
$$
(L_2-\chi_{2,2}^2\mathcal{E}+\chi_{1,2}^1\mathcal{E})\big(L_2+R_2-(\chi_{1,2}^1+\chi_{2,1}^1)\mathcal{E}\big)+(\chi_{2,2}^2-\chi_{1,2}^1-\chi_{2,1}^1)(R_2-\chi_{1,2}^1\mathcal{E})=0,
$$
i.e. $\mu(u_t,w_t)\in\langle u_t,w_t\rangle$. Analogously, if the required conditions are satisfied, then $\mu(w_t,u_t)\in\langle u_t,w_t\rangle$. We have also
\begin{multline*}
\mu(w_t,w_t)-\chi_{2,2}^2\alpha_t w_t=\alpha_t\big((L_2+R_2)M_t-\chi_{2,2}^2M_t\big)v\\
=\Big(\big((\chi_{1,2}^2+\chi_{2,1}^2-\chi_{1,1}^1)\mathcal{E}-S\big)S+(\chi_{1,2}^1-\chi_{2,1}^1-\chi_{2,2}^2)\mathcal{E}\Big)S\Big(\mathcal{E}+t\big((\chi_{1,2}^2+\chi_{2,1}^2-\chi_{1,1}^1)\mathcal{E}-S\big)\Big)=0
\end{multline*}
if the required conditions are satisfied.

It remains to prove that the subalgebra generated by $u=e_2+v$ has dimension not greater then $2$ if the equalities stated in the proposition are satisfied. Let $w=\mu(u,u)-\chi_{2,2}^2u=(L_2+R_2-\chi_{2,2}^2\mathcal{E})v$. Since $\mu(w,w)=0$, it remains to check that $\mu(u,w),\mu(w,u)\in\langle u,w\rangle$. Note that $w\in\langle e_3,\dots,e_n\rangle$ and $Sw=0$. Then we have
$$\mu(u,w)=L_2w=(\chi_{2,1}^2\mathcal{E}-R_1)Sw+\chi_{2,1}^1w=\chi_{2,1}^1w\mbox{ and }\mu(w,u)=R_2w=(\chi_{1,2}^2\mathcal{E}-L_1)Sw+\chi_{1,2}^1w=\chi_{1,2}^1w.$$
\end{Proof}

\begin{Remark} It follows from the just proved proposition that the classification of all trivial singular extension of a $2$-dimensional $1$-generated algebra $C$ is almost equivalent to the classification of pairs of matrices $(L,M)$ such that the minimal polynomial of $M$ divides some cubic polynomial defined by $C$.
In our problem we have more equivalences than in the classical problem of classification of matrices, but still not enough for making the problem solvable. By this reason, we do not give a general description of such algebras as in the case of  trivial singular extensions of $1$-dimensional algebras and  restrict our investigation by some estimations of levels.
\end{Remark}

Let us firstly consider the structures ${\bf k}^{n-2}\rtimes{\bf A}_3\in\mathcal{A}_n$ and ${\bf k}^2\rtimes_{\alpha,\beta}{\bf A}_3$ for $(\alpha,\beta)\in K_{1,1}^*$. Note that, for ${\bf k}=\mathbb{C}$, the algebras ${\bf k}^2\rtimes_{\alpha,\beta}{\bf A}_3$ ($(\alpha,\beta)\in K_{1,1}^*$) form the same set as the algebras $\mathfrak{N}_8$ and $\mathfrak{N}_9(\alpha)$ ($\alpha\in\mathbb{C}$) from \cite{kppv}.

\begin{Lem}\label{A3ext} $lev\left({\bf k}^{n-2}\rtimes{\bf A}_3\right)=3$ for any $n\ge 3$ and $lev\left({\bf k}^2\rtimes_{\alpha,\beta}{\bf A}_3\oplus {\bf k}^{n-4}\right)=3$ for any $(\alpha,\beta)\in K_{1,1}^*$ and $n\ge 4$.
\end{Lem}
\begin{Proof} Firstly, note that ${\bf k}^{n-2}\rtimes{\bf A}_3\xrightarrow{\left(e_1,\frac{1}{t}e_2+e_3,e_3,\dots,e_n\right)} T_0^{2,\overline{0,1}}$ and ${\bf k}^2\rtimes_{\alpha,\beta}{\bf A}_3\xrightarrow{(te_1,e_3,te_4te_2-e_4)} {\bf F}^{\alpha,\beta}\oplus {\bf k}$, and hence $lev\left({\bf k}^{n-2}\rtimes{\bf A}_3\right)\ge 3$ and $lev\left({\bf k}^2\rtimes_{\alpha,\beta}{\bf A}_3\right)\ge 3$.

Let us consider the set $$\mathcal{R}=\left\{\mu\in\mathcal{A}_n\left|\begin{array}{c}\mu(e_1,e_1)=\mu_{1,1}^ne_n,\,\mu(e_i,e_j)=0\mbox{ for $2\le i,j\le n$},\,\mu(e_1,e_n)=\mu(e_n,e_1)=0,\\
\mu(e_1,e_i)=\mu_{1,2}^2e_i+\mu_{1,i}^ne_n,\,\mu(e_i,e_1)=-\mu_{1,2}^2e_i-\mu_{1,i}^ne_n\mbox{ for $2\le i\le n-1$}\end{array}\right.\right\}.$$
It is easy to see that $\mathcal{R}$ is a closed subset of $\mathcal{A}_n$ invariant under lower triangular transformations of the basis $e_1,\dots,e_n$. It is also not difficult to see that ${\bf k}^n$, ${\bf A}_3\oplus {\bf k}^{n-2}$, ${\bf n}_3\oplus {\bf k}^{n-3}$, $T_0^{2,\overline{0,1}}$, ${\bf F}^{1,-1}\oplus {\bf k}^{n-3}$, and ${\bf k}^{n-2}\rtimes{\bf A}_3$ are all the structures, whose orbits intersect $\mathcal{R}$. Thus, $lev\left({\bf k}^{n-2}\rtimes{\bf A}_3\right)=3$.

Analogously, for any $(\alpha,\beta)\in K_{1,1}^*$, the set
$$\left\{\mu\in\mathcal{A}_n\left|\begin{array}{c} \mu(e_i,e_j)=0\,\mbox{if either $3\le i\le n$ or $3\le j\le n$},\,\mu(e_2,e_2)=0,\\
\mu_{1,2}^i=\mu_{2,1}^i=0\mbox{ for $1\le i\le n-1$},\,\alpha\mu_{2,1}^n=\beta\mu_{1,2}^n,\,\mu_{1,1}^i=0\mbox{ for $1\le i\le n-2$}\end{array}\right.\right\}$$
is a closed subset of $\mathcal{A}_n$ invariant under lower triangular transformations of the basis $e_1,\dots,e_n$. All the structures, whose orbits have a nontrivial intersection with this set, are ${\bf k}^n$, ${\bf A}_3\oplus {\bf k}^{n-2}$, ${\bf F}^{\alpha,\beta}\oplus {\bf k}^{n-3}$, and ${\bf k}^2\rtimes_{\alpha,\beta}{\bf A}_3\oplus {\bf k}^{n-4}$ if $\alpha+\beta\not=0$. If $\alpha+\beta=0$, then the set of structures with a nontrivial orbit intersection with the set under consideration contains additionally ${\bf n}_3$. In any case, we have $lev\left({\bf k}^2\rtimes_{\alpha,\beta}{\bf A}_3\oplus {\bf k}^{n-4}\right)=3$.
\end{Proof}

The next corollary is one of the results stated in \cite{khud15}.

\begin{corollary}\label{Flev} $lev\left({\bf k}^{n-3}\oplus {\bf F}^{\alpha,\beta}\right)=2$ for any $(\alpha,\beta)\in K_2^*$.
\end{corollary}
\begin{Proof} Follows from the proof of Lemma \ref{A3ext}.
\end{Proof}

\begin{Lem}\label{nondiag} Suppose that $\mu={\bf k}^{n-2}\rtimes_{(L_1,R_1,L_2,R_2)}\chi$ for some $\chi\in\mathcal{\tilde A}_2$ and $L_1,R_1,L_2,R_2\in M_{n-2}({\bf k})$.
Then either $\overline{O(\mu)}$ contains ${\bf k}^2\rtimes_{\alpha,\beta}{\bf A}_3\oplus {\bf k}^{n-4}$ for some  $(\alpha,\beta)\in K_{1,1}^*$ or $L_1=\alpha \mathcal{E}$ and $R_1=\beta \mathcal{E}$ for some $\alpha,\beta\in{\bf k}$.
\end{Lem}
\begin{Proof} If $L_1$ or $R_1$ is not diagonal, then $\mu(e_1,v)$ or $\mu(v,e_1)$ is linearly independent with $v$ for some $v\in\langle e_3,\dots,e_n\rangle$.
Thus, we may assume that $(\mu_{1,3}^4,\mu_{3,1}^4)\not=(0,0)$. Then $\mu\xrightarrow{(t^2e_1, t^4e_2, e_3, t^2e_4, te_5,\dots,te_n)}{\bf k}^2\rtimes_{\mu_{1,3}^4,\mu_{3,1}^4}{\bf A}_3\oplus {\bf k}^{n-4}.$
\end{Proof}

\begin{corollary}\label{notA3} If $A$ is a trivial singular extensions of ${\bf A}_3$ with $G(A)=2$, then $lev(A)\not=2$.
\end{corollary}
\begin{Proof} Suppose that the algebra $A$ is represented by ${\bf k}^{n-2}\rtimes_{(L_1,R_1,L_2,R_2)}{\bf A}_3$ for some matrices $L_1,R_1,L_2,R_2\in M_{n-2}({\bf k})$. If there are no $\alpha,\beta\in{\bf k}$ such that $L_1=\alpha \mathcal{E}$ and $R_1=\beta \mathcal{E}$, then $lev(A)\ge 3$ by Lemmas \ref{A3ext} and \ref{nondiag}. If $L_1=\alpha \mathcal{E}$ and $R_1=\beta \mathcal{E}$, then $(\alpha+\beta)^3=0$ by Proposition \ref{G2eq}. We have $L_2=R_2=0$ by the same proposition. Thus, $A$ can be represented by ${\bf A}_3\oplus{\bf k}^{n-2}$ if $\alpha=0$ and by ${\bf k}^{n-2}\rtimes{\bf A}_3$ if $\alpha\not=0$, and hence $lev(A)\in\{1,3\}$.
\end{Proof}

It follows from the results of \cite{kpv17} that the algebras presented in Table 4 with the exception of ${\bf E}_4$ are exactly all the $2$-dimensional $1$-generated algebra structures. In this table  we unite the series ${\bf E}_1$, ${\bf E}_2$ and ${\bf E}_3$ of the paper \cite{kpv17} in a one series called ${\bf E}_1$. We omit also the conditions required for the uniqueness modulo isomorphism, i.e. we allow some of the structures from Table 4 to represent the same algebra. Nevertheless, the structures that we really will consider, namely, ${\bf A}_1^{\alpha}$, ${\bf A}_2$, ${\bf B}_2^{\alpha}$ and ${\bf D}_2^{\alpha,\beta}$, where $\alpha,\beta\in{\bf k}$, $\alpha+\beta\not=1$, represent pairwise nonisomorphic algebras.
Note that some structures in Table 4 do not satisfy the conditions of Proposition \ref{G2eq}. To apply this proposition we have to apply firstly some linear transformation to the basis of $V$. On the other hand, any trivial singular extension is still determined by a structure of  $2$-dimensional $1$-generated algebra and four $(n-2)\times (n-2)$ matrices. We are going to classify trivial singular extension of $2$-dimensional $1$-generated algebra of the level $2$. By this reason, the next lemma allows to exclude $2$-dimensional algebras of the level $3$ from our consideration. Its proof is a direct calculation that we leave for the reader.

\begin{Lem}\label{3levext}
$$
\begin{aligned}
&{\bf k}^{n-2}\rtimes_{(L_1,R_1,L_2,R_2)}{\bf A}_4^{\alpha}&&\xrightarrow{(te_1-e_2,t^2e_2,e_3,\dots,e_n)}&&{\bf k}^{n-2}\rtimes_{(-L_2,-R_2,0,0)}{\bf A}_2,\\
&{\bf k}^{n-2}\rtimes_{(L_1,R_1,L_2,R_2)}{\bf B}_1^{\alpha}&&\xrightarrow{(e_1+te_2,-t^2e_2,e_3,\dots,e_n)}&&{\bf k}^{n-2}\rtimes_{(L_1,R_1,0,0)}{\bf A}_2,\\
&{\bf k}^{n-2}\rtimes_{(L_1,R_1,L_2,R_2)}{\bf C}^{\alpha,\beta}&&\xrightarrow{(te_1+e_2,t^2e_2,e_3,\dots,e_n)}&&{\bf k}^{n-2}\rtimes_{(L_2,R_2,0,0)}{\bf A}_1^{\alpha},\\
&{\bf k}^{n-2}\rtimes_{(L_1,R_1,L_2,R_2)}{\bf D}_1^{\alpha,\beta}&&\xrightarrow{(e_1,te_2,e_3,\dots,e_n)}&&{\bf k}^{n-2}\rtimes_{(L_1,R_1,0,0)}{\bf D}_2^{\beta,-\beta},\\
&{\bf k}^{n-2}\rtimes_{(L_1,R_1,L_2,R_2)}{\bf D}_3^{\alpha,\beta}&&\xrightarrow{(e_1,te_2,e_3,\dots,e_n)}&&{\bf k}^{n-2}\rtimes_{(L_1,R_1,0,0)}{\bf D}_2^{\alpha,\beta},\\
&{\bf k}^{n-2}\rtimes_{(L_1,R_1,L_2,R_2)}{\bf E}_1^{\alpha,\beta,\gamma,\delta}&&\xrightarrow{(e_1,te_2,e_3,\dots,e_n)}&&{\bf k}^{n-2}\rtimes_{(L_1,R_1,0,0)}{\bf D}_2^{\beta,\delta}.
\end{aligned}
$$
\end{Lem}

\begin{Remark} The parametrized bases in Lemma \ref{3levext} are taken from the proofs of \cite[Lemma 11, Theorem 13]{kpv17}. Note that one can use the same proofs to obtain in analogous way some additional degenerations for ${\bf B}_1^{\alpha}$, ${\bf D}_1^{\alpha,\beta}$, ${\bf D}_3^{\alpha,\beta}$ and ${\bf E}_1^{\alpha,\beta,\gamma,\delta}$.
\end{Remark}

Thus, it is enough to consider the trivial singular extensions of the structures ${\bf A}_1^{\alpha}$, ${\bf A}_2$, ${\bf B}_2^{\alpha}$ and ${\bf D}_2^{\alpha,\beta}$ for $\alpha,\beta\in{\bf k}$, $\alpha+\beta\not=1$.
Moreover, due to Lemmas \ref{A3ext} and \ref{nondiag}, in each case we may assume that $L_1=\alpha \mathcal{E}$ and $R_1=\beta \mathcal{E}$ for some $\alpha, \beta\in{\bf k}$.

\begin{Lem}\label{A1} Suppose that $\mu={\bf k}^{n-2}\rtimes_{(\epsilon \mathcal{E},\kappa \mathcal{E},L,R)}{\bf A}_1^{\alpha}$ for some $\alpha,\epsilon,\kappa\in {\bf k}$ and $L,R\in M_{n-2}({\bf k})$. If $G(\mu)=2$, then $\mu={\bf k}^{n-2}\rtimes_{\epsilon} {\bf A}_1^{\alpha}$. Moreover, $lev\big({\bf k}^{n-2}\rtimes_{\alpha} {\bf A}_1^{\alpha}\big)=2$ and $lev\big({\bf k}^{n-2}\rtimes_{\epsilon} {\bf A}_1^{\alpha}\big)\ge 3$ for $\epsilon\not=\alpha$.
\end{Lem}
\begin{Proof} It follows from Proposition \ref{G2eq} that $\epsilon+\kappa-1=0$ and $L=R=0$. Thus, $\mu={\bf k}^{n-2}\rtimes_{\epsilon} {\bf A}_1^{\alpha}$.

Let us now consider the set
$$\mathcal{R}=\left\{\chi\in\mathcal{A}_n\left|\begin{array}{c}\chi(e_1,e_1)=\chi_{1,1}^1e_1+\chi_{1,1}^ne_n,\,\chi(e_i,e_j)=0\mbox{ for $2\le i,j\le n$},\\
\chi(e_1,e_i)=\alpha\chi_{1,1}^1e_i,\chi(e_i,e_1)=(1-\alpha)\chi_{1,1}^1e_i\mbox{ for $2\le i\le n$}\end{array}\right.\right\}.$$
It is easy to see that $\mathcal{R}$ is a closed subset of $\mathcal{A}_n$ invariant under lower triangular transformations of the basis $e_1,\dots,e_n$. It is also not difficult to see that the set of orbits intersecting $\mathcal{R}$ is formed by the orbits of the structures ${\bf k}^n$, ${\bf A}_3\oplus {\bf k}^{n-2}$, $\nu^{\alpha}$, and ${\bf k}^{n-2}\rtimes_{\alpha} {\bf A}_1^{\alpha}$. Thus, $lev\big({\bf k}^{n-2}\rtimes_{\alpha} {\bf A}_1^{\alpha}\big)=2$. 
The degeneration ${\bf k}^{n-2}\rtimes_{\epsilon} {\bf A}_1^{\alpha}\xrightarrow{\left(e_1,\frac{1}{t}e_2+\frac{1}{\epsilon-\alpha}e_3,e_3,\dots,e_n\right)} T_1^{2,\overline{\alpha,\epsilon}}$ for $\epsilon\not=\alpha$ finishes the proof.
\end{Proof}

\begin{Lem}\label{A2} Suppose that $\mu={\bf k}^{n-2}\rtimes_{(\epsilon \mathcal{E},\kappa \mathcal{E},L,R)}{\bf A}_2$ for some $\epsilon,\kappa\in {\bf k}$ and $L,R\in M_{n-2}({\bf k})$. If $G(\mu)=2$, then $\mu={\bf k}^{n-2}\rtimes_{\epsilon} {\bf A}_2$. Moreover, $lev({\bf k}^{n-2}\rtimes_{1} {\bf A}_2)=2$ and $lev({\bf k}^{n-2}\rtimes_{\epsilon} {\bf A}_2)\ge 3$ for $\epsilon\not=1$.
\end{Lem}
\begin{Proof} It follows from Proposition \ref{G2eq} that $\epsilon+\kappa=0$ and $L=R=0$. Thus, $\mu={\bf k}^{n-2}\rtimes_{\epsilon} {\bf A}_2$.

Let us now consider the set
$$\mathcal{R}=\left\{\chi\in\mathcal{A}_n\left|\begin{array}{c}\chi(e_1,e_1)=\chi_{1,1}^ne_n,\,\chi(e_i,e_j)=0\mbox{ for $2\le i,j\le n$},\\
\chi(e_1,e_i)=\chi_{1,n}^ne_i,\,\chi(e_i,e_1)=-\chi_{1,n}^ne_i\mbox{ for $2\le i\le n$}\end{array}\right.\right\}.$$
It is easy to see that $\mathcal{R}$ is a closed subset of $\mathcal{A}_n$ invariant under lower triangular transformations of the basis $e_1,\dots,e_n$. It is also not difficult to see that the set of orbits intersecting $\mathcal{R}$ is formed by the orbits of the structures ${\bf k}^n$, ${\bf A}_3\oplus {\bf k}^{n-2}$, ${\bf p}^-$, and ${\bf k}^{n-2}\rtimes_{1} {\bf A}_2$. Thus, $lev({\bf k}^{n-2}\rtimes_{1} {\bf A}_2)=2$.
The degeneration ${\bf k}^{n-2}\rtimes_{\epsilon} {\bf A}_2\xrightarrow{\left(e_1,\frac{1}{t}e_2+\frac{1}{\epsilon-1}e_3,e_3,\dots,e_n\right)} T_0^{2,\overline{1,\epsilon}}$ for $\epsilon\not=1$ finishes the proof.
\end{Proof}

\begin{Lem}\label{B2} Suppose that $G(\mu)=2$ and $\mu={\bf k}^{n-2}\rtimes_{(\epsilon \mathcal{E}-L,\kappa \mathcal{E}-R,L,R)}{\bf B}_2^{\alpha}$ for some $\alpha,\epsilon,\kappa\in {\bf k}$ and $L,R\in M_{n-2}({\bf k})$.
 If $\epsilon=\alpha$, then either $\mu={\bf k}^{n-2}\rtimes_{\alpha} {\bf B}_2^{\alpha}$ or $\mu={\bf k}^{n-2}\rtimes_0^t {\bf B}_2^{\alpha}$. If $\epsilon\not=\alpha$, then either $\mu={\bf k}^{n-2}\rtimes_{\epsilon} {\bf B}_2^{\alpha}$ or $\mu\in O({\bf k}^{n-2}\rtimes_1^t {\bf B}_2^{\alpha})$.
Moreover, $lev\big({\bf k}^{n-2}\rtimes_0^t {\bf B}_2^{\alpha}\big)=2$, $lev\big({\bf k}^{n-2}\rtimes_{\epsilon} {\bf B}_2^{\alpha}\big)\ge 3$ and $lev\big({\bf k}^{n-2}\rtimes_1^t {\bf B}_2^{\alpha}\big)\ge 4$.
\end{Lem}
\begin{Proof} Let us consider the basis $e_1+e_2, e_2, e_3,\dots, e_n$ of $V$. Note that the structure constants $\tilde\mu_{i,j}^k$ ($1\le i,j,k\le n$) of $\mu$ in this basis satisfy the equalities $\tilde\mu_{2,2}^2=\tilde\mu_{i,j}^1=0$ for $1\le i,j\le 2$, $\tilde\mu_{1,1}^2=1$, $\tilde\mu_{1,2}^2=\alpha$, and $\tilde\mu_{2,1}^2=1-\alpha$. Thus, we can apply Proposition \ref{G2eq} and get that either $\epsilon+\kappa=0$ or  $\epsilon+\kappa=1$. In the first case we get $L=R=0$, i.e. $\mu={\bf k}^{n-2}\rtimes_{\epsilon} {\bf B}_2^{\alpha}$.
In the second case we get $L=(\epsilon-\alpha)\mathcal{E}$, $R=(\alpha-\epsilon)\mathcal{E}$. If $\epsilon=\alpha$, then $\mu={\bf k}^{n-2}\rtimes_0^t {\bf B}_2^{\alpha}$. If $\epsilon\not=\alpha$, then, rescaling $e_2$ by $\frac{1}{\epsilon-\alpha}$, we get $\mu\in O({\bf k}^{n-2}\rtimes_1^t {\bf B}_2^{\alpha})$.

Let us now consider the set
$$\mathcal{R}=\left\{\chi\in\mathcal{A}_n\left|\begin{array}{c}\chi_{1,1}^1=0,\,\chi(e_i, e_j)=0\mbox{ for $2\le i,j\le n$},\\\chi(e_1,e_i)=\alpha(\chi_{1,2}^2+\chi_{2,1}^2)e_i,\,\chi(e_i,e_1)=(1-\alpha)(\chi_{1,2}^2+\chi_{2,1}^2)e_i\mbox{ for $2\le i\le n$}\end{array}\right.\right\}.$$
Direct verification shows that $\mathcal{R}$ is a closed subset of $\mathcal{A}_n$ invariant under lower triangular transformations of the basis $e_1,\dots,e_n$.
The set of orbits intersecting $\mathcal{R}$ is formed by the orbits of the structures ${\bf k}^n$, ${\bf A}_3\oplus {\bf k}^{n-2}$ and ${\bf k}^{n-2}\rtimes_{0}^t {\bf B}_2^{\alpha}$. Thus, $lev\big({\bf k}^{n-2}\rtimes_0^t {\bf B}_2^{\alpha}\big)=2$.

Since
$$
    {\bf k}^{n-2}\rtimes_{\epsilon} {\bf B}_2^{\alpha}\xrightarrow{(te_1+e_2, e_2+e_3, te_2,e_4,\dots,e_n)} {\bf F}^{\alpha-\epsilon,1-\alpha+\epsilon}\oplus{\bf k}^{n-3}\mbox{ and }
    {\bf k}^{n-2}\rtimes_1^t {\bf B}_2^{\alpha}\xrightarrow{(te_1+e_2,te_2,e_3,\dots,e_n)} {\bf k}^{n-2}\rtimes {\bf A}_3
$$
we have $lev\big({\bf k}^{n-2}\rtimes_{\epsilon} {\bf B}_2^{\alpha}\big)\ge 3$ for any $\epsilon\in{\bf k}$ and $lev\big({\bf k}^{n-2}\rtimes_1^t {\bf B}_2^{\alpha}\big)\ge 4$.
\end{Proof}

\begin{Lem}\label{D2} Suppose that $G(\mu)=2$ and $\mu={\bf k}^{n-2}\rtimes_{(\epsilon \mathcal{E}-L,\kappa \mathcal{E}-R,L,R)}{\bf D}_2^{\alpha,\beta}$ for some $\alpha,\beta,\epsilon,\kappa\in {\bf k}$ ($\alpha+\beta\not=1$) and $L,R\in M_{n-2}({\bf k})$.
 If $\epsilon=\alpha$, then either $\mu= {\bf k}^{n-2}\rtimes_{\alpha} {\bf D}_2^{\alpha,\beta}$ or $\mu={\bf k}^{n-2}\rtimes_0^t {\bf D}_2^{\alpha,\beta}$.
If $\epsilon\not=\alpha$, then either $\mu={\bf k}^{n-2}\rtimes_{\epsilon} {\bf D}_2^{\alpha,\beta}$ or $\mu\in O({\bf k}^{n-2}\rtimes_1^t {\bf D}_2^{\alpha,\beta})$.
Moreover, $lev\big({\bf k}^{n-2}\rtimes_0^t {\bf D}_2^{\alpha,\beta}\big)=2$, $lev\big({\bf k}^{n-2}\rtimes_{\epsilon} {\bf D}_2^{\alpha,\beta}\big)\ge 3$ and $lev\big({\bf k}^{n-2}\rtimes_1^t {\bf D}_2^{\alpha,\beta}\big)\ge 4$.
\end{Lem}
\begin{Proof} Let us consider the basis $e_1+e_2, (\alpha+\beta-1)e_2, e_3,\dots, e_n$ of $V$. Note that the structure constants $\tilde\mu_{i,j}^k$ ($1\le i,j,k\le n$) of $\mu$ in this basis satisfy the equalities $\tilde\mu_{2,2}^2=\tilde\mu_{1,2}^1=\tilde\mu_{2,1}^1=0$, $\tilde\mu_{1,1}^1=\tilde\mu_{1,1}^2=1$, $\tilde\mu_{1,2}^2=\alpha$, and $\tilde\mu_{2,1}^2=\beta$. Thus, we can apply Proposition \ref{G2eq} and get that either $\epsilon+\kappa=1$ or  $\epsilon+\kappa=\alpha+\beta$.
In the first case we get $L=R=0$, i.e. $\mu={\bf k}^{n-2}\rtimes_{\epsilon} {\bf D}_2^{\alpha,\beta}$.
In the second case we get $L=(\epsilon-\alpha)\mathcal{E}$, $R=(\alpha-\epsilon)\mathcal{E}$. If $\epsilon=\alpha$, then $\mu={\bf k}^{n-2}\rtimes_0^t {\bf D}_2^{\alpha,\beta}$. If $\epsilon\not=\alpha$, then, rescaling $e_2$ by $\frac{1}{\epsilon-\alpha}$, one can see that $\mu\in O({\bf k}^{n-2}\rtimes_1^t {\bf D}_2^{\alpha,\beta})$.

Let us now consider the set
$$\mathcal{R}=\left\{\chi\in\mathcal{A}_n\left|\begin{array}{c}\chi(e_i, e_j)=0\mbox{ for $2\le i,j\le n$},\chi(e_1,e_i)=\alpha\chi_{1,1}^1e_i,\,\chi(e_i,e_1)=\beta\chi_{1,1}^1e_i\mbox{ for $2\le i\le n$}\end{array}\right.\right\}.$$
Direct verification shows that $\mathcal{R}$ is a closed subset of $\mathcal{A}_n$ invariant under lower triangular transformations of the basis $e_1,\dots,e_n$.
The set of orbits intersecting $\mathcal{R}$ is formed by the orbits of the structures ${\bf k}^n$, ${\bf A}_3\oplus {\bf k}^{n-2}$ and ${\bf k}^{n-2}\rtimes_{0}^t {\bf D}_2^{\alpha,\beta}$. Thus, $lev\big({\bf k}^{n-2}\rtimes_0^t {\bf D}_2^{\alpha,\beta}\big)=2$.

Since
\begin{multline*}
    {\bf k}^{n-2}\rtimes_{\epsilon} {\bf D}_2^{\alpha,\beta}\xrightarrow{\left(te_1+e_2, (\alpha+\beta-1)(e_2+e_3),(\alpha+\beta-1)te_2,e_4,\dots,e_n\right)} {\bf F}^{\alpha-\epsilon,\beta+\epsilon-1}\oplus{\bf k}^{n-3}\mbox{ and }\\
    {\bf k}^{n-2}\rtimes_1^t {\bf D}_2^{\alpha,\beta}\xrightarrow{\left(te_1+e_2, (\alpha+\beta-1)te_2,e_3,\dots,e_n\right)} {\bf k}^{n-2}\rtimes {\bf A}_3
\end{multline*}
we have $lev\big({\bf k}^{n-2}\rtimes_{\epsilon} {\bf D}_2^{\alpha,\beta}\big)=3$ for any $\epsilon\in{\bf k}$ and $lev\big({\bf k}^{n-2}\rtimes_1^t {\bf D}_2^{\alpha,\beta}\big)\ge 4$.
\end{Proof}

\subsection{Algebras with ${\bf A}_3$-ideal}
This subsection is devoted to the case where IW contraction with respect to a $2$-dimensional $1$-generated algebra gives an algebra of level $1$, namely, to algebras that have an ideal isomorphic to ${\bf A}_3$ as an algebra. The first important example of such an algebra is the algebra of a bilinear form.

\begin{Def}{\rm
An {\it algebra of a bilinear form} is an algebra that can be represented by a structure $\mu\in\mathcal{A}_n$ such that $\mu_{i,j}^k$ have nonzero values only for $1\le i,j\le n-1$ and $k=n$. If at the same time $\mu(v,v)=0$ for all $v\in V$, then we call the corresponding algebra an {\it algebra of an antisymmetric bilinear form}. In the opposite case we call it an {\it algebra of a nonantisymmetric bilinear form}.
}
\end{Def}

It follows from the classification of antisymmetric bilinear forms that any algebra of an antisymmetric bilinear form is either trivial or can be represented by $\eta_m\oplus {\bf k}^{n-2m-1}$ for some $1\le m< \frac{n}{2}$.
It is easy to see also that an algebra of a nonantisymmetric bilinear form contains an ideal isomorphic to ${\bf A}_3$ as an algebra. We give here an estimation of the level of an algebra of a bilinear. Let us recall that, for an algebra $A$, the {\it annihilator} of $A$ is the ideal $Ann(A)=\{a\in A\mid ab=ba=0\mbox{ for all $b\in A$}\}$.

\begin{proposition}\label{levelform}
Let $A$ be an algebra of a nonantisymmetric bilinear form. Then $lev(A)\ge n-dim\,Ann(A)$.
\end{proposition}
\begin{Proof} Let us prove the assertion by induction on $n-dim\,Ann(A)$. If $n-dim\,Ann(A)=1$, then it is easy to see that $A$ can be represented by ${\bf A}_3\oplus{\bf k}^{n-2}$, and hence $lev(A)=1$. Let now put $m=n-dim\,Ann(A)>1$. We may assume that $A$ is represented by a structure $\mu\in\mathcal{A}_n$ such that  $\mu(e_m,e_m)=1$ and $\mu_{i,j}^k$ takes nonzero values only for $1\le i,j\le m$ and $k=n$. Let us transform the elements $e_1,\dots, e_m$ in the following way. Firstly, let us replace $e_i$ by $e_i-\mu_{i,m}^ne_m$ for $1\le i\le m-1$. After this replacement we may assume that $\mu_{i,m}^n=0$ for $1\le i\le m-1$. Now, if $\mu_{m,i}^n\not=0$ for some $1\le i\le m-1$, then, after an obvious manipulation, we may assume that $\mu_{m,m-1}^n=1$. Now, replacing $e_i$ by $e_i-\mu_{m,i}^ne_{m-1}$, we may assume that $\mu_{m,i}^n=0$ for $1\le i\le m-2$. Let us now choose $\alpha\in{\bf k}$ such that $\mu(e_{m-1}+\alpha e_m,e_{m-1}+\alpha e_m)\not=0$ and replace $e_{m-1}$ by $e_{m-1}+\alpha e_m$ for such a scalar $\alpha$. Thus, we may assume that $\mu(e_{m-1},e_{m-1})\not=0$ and $\mu(e_i,e_m)=\mu(e_i,e_m)=0$ for any $1\le i\le m-2$. Let us now consider two cases.
\begin{enumerate}
\item For any nonzero $u\in\langle e_1,\dots, e_{m-1}\rangle$ there exists $1\le i\le m-1$ such that either $\mu(e_i,u)\not=0$ or $\mu(u,e_i)\not=0$. Let us consider the degeneration $A\to B$ defined by the parametrized basis given by the equalities $E_i^t=e_i$ for $1\le i\le n$, $i\not=m$ and $E_m^t=te_m$. It is easy to see that $B$ is an algebra of a bilinear form such that $a^2\not=0$ for some $a\in B$ and $n-dim\,Ann(B)=m-1$. Then by induction hypothesis we have $lev(B)\ge m-1$, and hence $lev(A)\ge m$.
\item There exists nonzero $u\in\langle e_1,\dots, e_{m-1}\rangle$ such that $\mu(e_i,u)=\mu(u,e_i)=0$ for $1\le i\le m-1$. If $u\in\langle e_1,\dots, e_{m-2}\rangle$, then $\mu(u,e_m)=\mu(e_m,u)=0$, and hence $dim\,Ann(A)>n-m$. Thus, $u\not\in\langle e_1,\dots, e_{m-2}\rangle$ and we can replace $e_{m-1}$ by $u$. Note that we still have $\mu(e_i,e_m)=\mu(e_m,e_i)=0$ for $1\le i\le m-2$ while the condition $\mu(e_{m-1},e_{m-1})\not=0$ may be not satisfied. On the other hand, now we have $\mu(e_i,e_{m-1})=\mu(e_{m-1},e_i)=0$ for $1\le i\le m-2$. Let us prove that for any nonzero $v\in\langle e_1,\dots, e_{m-2},e_m\rangle$ there exists $i\in\{1,\dots,m-2\}\cup\{m\}$ such that either $\mu(e_i,v)\not=0$ or $\mu(v,e_i)\not=0$. Really, if $v\in\langle e_1,\dots, e_{m-2},e_m\rangle$ does not satisfy the required condition, then it follows from $\mu(v,e_m)=0$ that $v\in\langle e_1,\dots, e_{m-2}\rangle$, and hence $\mu(e_{m-1},v)=\mu(v,e_{m-1})=0$. This again contradicts to the fact that $dim\,Ann(A)=n-m$. Let us consider the degeneration $A\to B$ defined by the parametrized basis given by the equalities $E_i^t=e_i$ for $1\le i\le n$, $i\not=m-1$ and $E_{m-1}^t=te_{m-1}$. Applying the induction hypothesis to $B$ we get $lev(A)\ge lev(B)+1\ge m$.
\end{enumerate}
\end{Proof}

Suppose that the algebra $A$ has an ideal isomorphic to ${\bf A}_3$ as an algebra. Let us represent $A$ by a structure $\mu\in\mathcal{A}_n$ such that $\mu(e_{n-1},e_{n-1})=e_n$, $\mu(e_{n-1},e_n)=\mu(e_{n},e_{n-1})=\mu(e_{n},e_n)=0$ and $\mu(e_i,\langle e_{n-1},e_n\rangle),\mu(\langle e_{n-1},e_n\rangle,e_i)\subset\langle e_{n-1},e_n\rangle$ for any $1\le i\le n-2$. Then there is a degeneration $A\to B$ corresponding to the parametrized basis defined by the equalities $E_i^t=te_i$ for $1\le i\le n-1$ and $E_n^t=t^2e_n$. It is easy to see that $B$ is an algebra of a nonantisymmetric bilinear form. We will call $B$ an {\it ${\bf A}_3$-bilinear form contraction} of $A$. Our next goal is to describe all the algebras with an ideal isomorphic to ${\bf A}_3$, whose ${\bf A}_3$-bilinear form contractions are of the level $1$. All of these algebras except ${\bf A}_3\oplus{\bf k}^{n-2}$ can be found in Table 6.

\begin{Lem}\label{nonform_class} Suppose that $A$ has an ideal isomorphic to ${\bf A}_3$ as an algebra. If all the ${\bf A}_3$-bilinear form contractions of $A$ can be represented by ${\bf A}_3\oplus{\bf k}^{n-2}$, then $A$ can be represented by a structure from the set
$$\{{\bf A}_3\oplus{\bf k}^{n-2},{\bf A}_3\rtimes{\bf p}^-,{\bf A}_3\rtimes \nu^{\alpha},({\bf A}_3\oplus{\bf k}^{n-4})\rtimes {\bf E}_4\}.$$
\end{Lem}
\begin{Proof} It is enough to consider the case $n\ge 3$. Let us represent $A$ by a structure $\mu$ that is described just before this proposition. Suppose that all the ${\bf A}_3$-bilinear form contraction of $A$ can be represented by ${\bf A}_3$. Let us now replace $e_i$ by $e_i-\mu_{i,n-1}^ne_{n-1}$ for all $1\le i\le n-2$. After this replacement we may assume that $\mu_{i,n-1}^n=0$ for $1\le i\le n-2$. It is easy to see that if $\mu_{i,j}^n\not=0$ or $\mu_{n-1,i}^n\not=0$ for some $1\le i,j\le n-2$, then the ${\bf A}_3$-bilinear form contraction of $A$ corresponding to the structure $\mu$ has the dimension of the annihilator less than $n-1$, i.e. cannot be represented by ${\bf A}_3\oplus{\bf k}^{n-2}$.
Thus, $\mu_{i,j}^n=\mu_{n-1,i}^n=0$ for all $1\le i,j\le n-2$

Let $1\le m\le n-2$ be some integer and $\tilde\mu_{i,j}^k$ ($1\le i,j,k\le n$) be the structure constants of $\mu$ in the basis $e_1,\dots,e_{m-1},e_m-e_n,e_{m+1},\dots,e_n$ of $V$. Then it is easy to see that $\tilde\mu_{i,j}^n=\mu_{i,j}^m$ for all $1\le i,j\le n-2$, $i,j\not=m$. As it was mentioned above, it follows that $\mu_{i,j}^k=0$ for all $1\le i,j,k\le n-2$ such that $i,j\not=k$. This condition has to be also satisfied after any nondegenerate linear transformation of the elements $e_1,\dots,e_{n-2}$. In other words, $ab\in\langle a,b\rangle+{\bf A}_3$ for any $a,b\in A$, where the ideal ${\bf A}_3$ corresponds to the subspace $\langle e_{n-1},e_n\rangle$ of $V$ after going from $A$ to $\mu$. In particular, $G(A/{\bf A}_3)=1$. Let us take some element $a\in A/{\bf A}_3$ and consider the IW contraction $A/{\bf A}_3\to B$ with respect to $(A/{\bf A}_3)(a)$. Due to the results of Section \ref{TnGen}, the algebra $B$ can be represented by some $T_r^M\in\mathcal{A}_{n-2}$, where $r\in\{0,1\}$ and $M\in M_{n-3}({\bf k})$. Since $bc$ has to lie in $\langle b,c\rangle$ for any $b,c\in B$, the matrix $M$ has to be diagonal in any basis. Thus, the algebra $A/{\bf A}_3$ satisfies the conditions of Corollary \ref{nu_ptype}, and hence either is trivial or can be represented by a structure from the set $\{{\bf p}^-,\nu^{\alpha},{\bf k}^{n-2}\rtimes {\bf E}_4\}_{\alpha\in {\bf k}}$.

Let now $\tilde\mu_{i,j}^k(\alpha)$ ($1\le i,j,k\le n$, $\alpha\in{\bf k}$) denote the structure constants of $\mu$ in the basis $$e_1+\big(\alpha(\mu_{1,n}^n-\mu_{1,n-1}^{n-1})+\alpha^2\mu_{1,n}^{n-1}\big)e_{n-1},\dots,e_{n-2}+\big(\alpha(\mu_{n-2,n}^n-\mu_{n-2,n-1}^{n-1})+\alpha^2\mu_{n-2,n}^{n-1}\big)e_{n-1},e_{n-1}-\alpha e_n,e_n.$$
Direct calculation shows that $\tilde\mu_{i,n-1}^n(\alpha)=0$ for $1\le i\le n-2$, and hence $\tilde\mu_{i,j}^n(\alpha)$ and $\tilde\mu_{n-1,i}^n(\alpha)$ have to be zero for all $1\le i,j\le n-2$ and all $\alpha\in{\bf k}$; in particular, $\tilde\mu_{i,i}^n(\alpha)$ has to be zero. On the other hand, we have
\begin{multline*}
\tilde\mu_{i,i}^n(\alpha)=\alpha\mu_{i,i}^{n-1}+\alpha^2(\mu_{i,n}^n-\mu_{i,n-1}^{n-1})(\mu_{i,n}^n+\mu_{n-1,i}^{n-1}-\mu_{i,i}^i)\\
+\alpha^3\big(\mu_{i,n}^{n-1}\mu_{n-1,i}^{n-1}+\mu_{i,n}^n\mu_{i,n}^{n-1}+\mu_{i,n}^{n-1}(\mu_{i,n}^n-\mu_{i,n-1}^{n-1})-\mu_{i,i}^i\mu_{i,n}^{n-1}\big)+\alpha^4\left(\mu_{i,n}^{n-1}\right)^2,\\
\tilde\mu_{n-1,i}^n(\alpha)=\alpha(\mu_{n-1,i}^{n-1}-\mu_{n,i}^n+\mu_{i,n}^n-\mu_{i,n-1}^{n-1})+\alpha^2(\mu_{i,n}^{n-1}-\mu_{n,i}^{n-1}).
\end{multline*}
We get that $\mu_{i,i}^{n-1}=0$ and $\mu_{i,n}^{n-1}=0$ for all $1\le i\le n-2$ from the first equality. Then it follows from the second equality that $\mu_{n,i}^{n-1}=0$ too. It is not difficult also to deduce that 
either $\mu_{i,n-1}^{n-1}=\mu_{i,n}^n$, $\mu_{n-1,i}^{n-1}=\mu_{n,i}^n$ or $\mu_{n-1,i}^{n-1}+\mu_{i,n}^n=\mu_{i,n-1}^{n-1}+\mu_{n,i}^n=\mu_{i,i}^i$.

 Now, for $1\le i,j\le n-2$, $i\not=j$, we have
$$
\tilde\mu_{i,j}^n(\alpha)=\alpha\mu_{i,j}^{n-1}+\alpha^2\big((\mu_{i,n}^n(\mu_{n-1,j}^{n-1}+\mu_{j,n}^n-\mu_{j,n-1}^{n-1})-\mu_{i,n-1}^{n-1}\mu_{n-1,j}^{n-1}-\mu_{i,j}^i(\mu_{i,n}^n-\mu_{i,n-1}^{n-1})-\mu_{i,j}^j(\mu_{j,n}^n-\mu_{j,n-1}^{n-1})\big).
$$
Thus, $\mu_{i,j}^{n-1}=0$ for all $1\le i,j\le n-2$.

As it was proved above, $A/{\bf A}_3$ either is trivial or can be represented by a structure from the set $\{{\bf p}^-,\nu^{\alpha},{\bf k}^{n-2}\rtimes {\bf E}_4\}_{\alpha\in {\bf k}}$. Let us consider all the cases separately.
\begin{enumerate}
\item Suppose that $A/{\bf A}_3$ is trivial. Let us consider firstly the case $n>3$. Let us prove that $\mu_{i,n}^n=\mu_{n,i}^n=0$ for any $1\le i\le n-2$. Let us choose some $1\le m\le n-2$ such that $m\not=i$ and consider the basis $e_1,\dots,e_{m-1},e_m+e_n,e_{m+1},\dots,e_n$ of $V$. If $\mu_{i,n}^n\not=0$ or $\mu_{n,i}^n\not=0$, then it is easy to see that the structure constants in this basis do not satisfy all the conditions obtained above since $\mu(e_i,e_m+e_n)=\mu_{i,n}^ne_n$ and $\mu(e_m+e_n,e_i)=\mu_{n,i}^ne_n$. Thus, $A$ can be represented by ${\bf A}_3\oplus{\bf k}^{n-2}$ if $n>3$. If $n=3$, then considering the basis $e_1+e_3,e_2,e_3$ we get that $\mu_{1,3}^3+\mu_{3,1}^3=0$. Then it is easy to see that $A$ can be represented either by the structure ${\bf A}_3\oplus{\bf k}$ or by the structure ${\bf A}_3\rtimes{\bf p}^-$.
\item Suppose that $A/{\bf A}_3$ can be represented by ${\bf p}^-$; in particular, $n\ge 4$. Then we may assume that $\mu(e_1,e_1)=0$, $\mu(e_1,e_i)=-\mu(e_i,e_1)=e_i$ and $\mu(e_i,e_j)=0$ for $2\le i,j\le n-2$. Let $\tilde\mu_{i,j}^k$ ($1\le i,j,k\le n$) be the structure constants of $\mu$ in the basis $e_1,e_2+e_n,e_3,\dots,e_n$ of $V$. It is easy to see that $\tilde\mu_{1,2}^n=\mu_{1,n}^n-1$ and $\tilde\mu_{2,1}^n=\mu_{n,1}^n+1$. Thus, $\mu_{1,n}^n=-\mu_{n,1}^n=1$, and hence $\mu_{1,n-1}^{n-1}=-\mu_{n-1,1}^{n-1}=1$. Analogously, considering the basis $e_1+e_n,e_2,e_3,\dots,e_n$ of $V$, we get $\mu(e_i,e_{n-1})=\mu(e_{n-1},e_i)=\mu(e_i,e_n)=\mu(e_n,e_i)=0$ for $2\le i\le n-2$. Thus, $\mu={\bf A}_3\rtimes{\bf p}^-$.
\item Suppose that $A/{\bf A}_3$ can be represented by $\nu^{\alpha}$ for some $\alpha\in{\bf k}$. Then we may assume that $\mu(e_1,e_1)=e_1$, $\mu(e_1,e_i)=\alpha e_i$, $\mu(e_i,e_1)=(1-\alpha)e_i$ and $\mu(e_i,e_j)=0$ for $2\le i,j\le n-2$. As in the previous case, considering the basis $e_1+e_n,e_2,e_3,\dots,e_n$ of $V$, we get $\mu(e_i,e_{n-1})=\mu(e_{n-1},e_i)=\mu(e_i,e_n)=\mu(e_n,e_i)=0$ for $2\le i\le n-2$ and $\mu_{1,n}^n+\mu_{n,1}^n=1$. If $n=3$, then in fact the algebra $A/{\bf A}_3$ does not depend on $\alpha$ and we can set $\alpha=\mu_{1,n}^n$. If $n>3$, then, considering the basis $e_1,e_2+e_n,e_3,\dots,e_n$ of $V$, we get $\mu_{1,n}^n=\alpha$. In any case, we have $\mu_{1,n}^n=\alpha$, $\mu_{n,1}^n=1-\alpha$ and, taking in account the equalities proved above, $\mu_{1,n-1}^{n-1}=\alpha$, $\mu_{n-1,1}^{n-1}=1-\alpha$. Thus, $\mu={\bf A}_3\rtimes \nu^{\alpha}$.
\item Suppose that $A/{\bf A}_3$ can be represented by ${\bf k}^{n-2}\rtimes {\bf E}_4$. Then, analogously to the previous case, we get that $\mu=({\bf A}_3\oplus{\bf k}^{n-4})\rtimes {\bf E}_4$.
\end{enumerate}
\end{Proof}

Now it is not very difficult to describe all the algebras of the level $2$ with an ideal isomorphic to ${\bf A}_3$. It is also possible to calculate the exact values of the levels of ${\bf A}_3\rtimes{\bf p}^-$, ${\bf A}_3\rtimes \nu^{\alpha}$ and $({\bf A}_3\oplus{\bf k}^{n-4})\rtimes {\bf E}_4$, but in this paper we will show only that they all have levels not less than $3$. This fact follows from the next technical lemma, whose proof we leave for the reader.

\begin{Lem}\label{nonform} For any $\alpha\in{\bf k}$ we have degenerations
$$
\begin{aligned}
&({\bf A}_3\oplus{\bf k}^{n-4})\rtimes {\bf E}_4&&\xrightarrow{(e_1-e_2,te_2,e_3,\dots,e_n)}&&{\bf A}_3\rtimes{\bf p}^-&&\xrightarrow{(e_1+e_{n-1},e_n,te_2,\dots,te_{n-1})}&&{\bf k}^{n-2}\rtimes_{1} {\bf A}_2,\\
&({\bf A}_3\oplus{\bf k}^{n-4})\rtimes {\bf E}_4&&\xrightarrow{(\alpha e_1+(1-\alpha)e_2,te_1-te_2,e_3,\dots,e_n)}&&{\bf A}_3\rtimes_{\alpha}\nu^{\alpha}&&\xrightarrow{(e_1+e_{n-1},e_n,te_2,\dots,te_{n-1})}&&{\bf k}^{n-2}\rtimes_{\alpha} {\bf A}_1^{\alpha}.
\end{aligned}
$$
\end{Lem}

\begin{proposition}\label{A3ideal} Suppose that $A$ has an ideal isomorphic to ${\bf A}_3$ as an algebra. Then it has level $2$ iff it can be represented by some structure from the set $\{{\bf F}^{\alpha,\beta}\oplus{\bf k}^{n-3}\}_{(\alpha,\beta)\in K_{2}^*}$.
\end{proposition}
\begin{Proof} By Lemmas \ref{nonform_class} and \ref{nonform}, any algebra with an ideal isomorphic to ${\bf A}_3$ of the level $2$ can be contracted to an algebra of a nonantisymmetric bilinear form with level greater or equal to $2$. Thus, $A$ can have the level $2$ only if it is an algebra of a nonantisymmetric bilinear form. Moreover, by Proposition \ref{levelform}, we have $dim\,Ann(A)\ge n-2$. It is easy to see that if $dim\,Ann(A)> n-2$, then $A$ either is trivial or can be represented by ${\bf A}_3\oplus{\bf k}^{n-2}$. On the other hand, any algebra of  a nonantisymmetric bilinear form with annihilator of codimension $2$ can be represented by the structure ${\bf F}^{\alpha,\beta}\oplus{\bf k}^{n-3}$ for some $(\alpha,\beta)\in K_{2}^*$, whose level is equal to $2$ by Corollary \ref{Flev}.
\end{Proof}

Note that, for ${\bf k}=\mathbb{C}$, the algebras ${\bf F}^{\alpha,\beta}\oplus {\bf k}$ ($(\alpha,\beta)\in K_2^*$) form the same set as the algebras $\mathfrak{N}_2^{\mathbb{C}}(\beta)$ ($\beta\in\mathbb{C}$) and $\mathfrak{N}_3^{\mathbb{C}}$ from \cite{kppv} and  the algebras ${\bf F}^{\alpha,\beta}\oplus {\bf k}^{n-3}$ ($(\alpha,\beta)\in K_2^*$) form the same set as the algebras $A_5(\alpha)$ ($\alpha\in\mathbb{C}\setminus\{-1\}$) and $A_6$ from \cite{khud15}.

\subsection{Classification of algebras of the level $2$}

In this subsection we apply the results of previous sections to get a classification of the algebras of the level $2$. As a corollary, in the end of this subsection we will give the same classification in some certain varieties. In particular, we will recover the results of \cite{khud15} and will generalize some results of \cite{gorb93}.
Thus, the main result of this subsection and one of the main results of the present paper is the next theorem.

\begin{Th}\label{ClassLev2} Let $A$ be an $n$-dimensional algebra.
\begin{enumerate}
\item If $n=2$, then $lev(A)=2$ iff $A$ can be represented by some structure from the set
$$\{{\bf A}_1^{\alpha},{\bf B}_2^{\alpha}\}_{\alpha\in{\bf k}}\cup\{{\bf D}_2^{\alpha,\beta}\}_{\alpha,\beta\in{\bf k},\alpha+\beta\not=1}\cup\{{\bf A}_2,{\bf E}_4\}.$$
\item If $n=3$, then $lev(A)=2$ iff $A$ can be represented by a structure from the set
\begin{multline*}
\{{\bf k}\rtimes_{\alpha}{\bf A}_1^{\alpha},{\bf k}\rtimes_0^t{\bf B}_2^{\alpha}\}_{\alpha\in{\bf k}}\cup\{{\bf k}\rtimes_0^t{\bf D}_2^{\alpha,\beta}\}_{\alpha,\beta\in{\bf k},\alpha+\beta\not=1}\cup\{{\bf F}^{\alpha,\beta}\}_{(\alpha,\beta)\in K_{2}^*}\\
\cup\left\{T_0^{2,\overline{\alpha,\beta}}\right\}_{(\alpha,\beta)\in K_2^*}\cup\left\{T_1^{2,\overline{\alpha,\beta}}\right\}_{(\alpha,\beta)\in K_2}\cup\{{\bf k}\rtimes_{1} {\bf A}_2,{\bf k}\rtimes {\bf E}_4\}.
\end{multline*}
\item If $n=4$, then $lev(A)=2$ iff $A$ can be represented by a structure from the set
\begin{multline*}
\{{\bf k}^2\rtimes_{\alpha}{\bf A}_1^{\alpha},{\bf k}^2\rtimes_0^t{\bf B}_2^{\alpha}\}_{\alpha\in{\bf k}}\cup\{{\bf k}^2\rtimes_0^t{\bf D}_2^{\alpha,\beta}\}_{\alpha,\beta\in{\bf k},\alpha+\beta\not=1}\cup\{{\bf F}^{\alpha,\beta}\oplus{\bf k}\}_{(\alpha,\beta)\in K_{2}^*}\\
\cup\left\{T_0^{2,\overline{\alpha,\beta}}\right\}_{(\alpha,\beta)\in K_{1,1}^*}\cup\left\{T_1^{2,\overline{\alpha,\beta}}\right\}_{(\alpha,\beta)\in K_{1,1}}\cup\{{\bf k}^2\rtimes_{1} {\bf A}_2,T_{0}^3,{\bf k}^2\rtimes {\bf E}_4\}.
\end{multline*}
\item If $n\ge 5$, then $lev(A)=2$ iff $A$ can be represented by a structure from the set
\begin{multline*}
\{{\bf k}^{n-2}\rtimes_{\alpha}{\bf A}_1^{\alpha},{\bf k}^{n-2}\rtimes_0^t{\bf B}_2^{\alpha}\}_{\alpha\in{\bf k}}\cup\{{\bf k}^{n-2}\rtimes_0^t{\bf D}_2^{\alpha,\beta}\}_{\alpha,\beta\in{\bf k},\alpha+\beta\not=1}\cup\{{\bf F}^{\alpha,\beta}\oplus{\bf k}^{n-3}\}_{(\alpha,\beta)\in K_{2}^*}\\
\cup\left\{T_0^{2,\overline{\alpha,\beta}}\right\}_{(\alpha,\beta)\in K_{1,1}^*}\cup\left\{T_1^{2,\overline{\alpha,\beta}}\right\}_{(\alpha,\beta)\in K_{1,1}}\cup\{{\bf k}^{n-2}\rtimes_{1} {\bf A}_2,T_{0}^{2,2},{\bf k}^{n-2}\rtimes {\bf E}_4,\eta_2\oplus{\bf k}^{n-5}\}.
\end{multline*}
\end{enumerate}
\end{Th}
\begin{Proof} It follows from Propositions \ref{GenType1}, \ref{A3ideal} and Lemmas \ref{A1}, \ref{A2}, \ref{B2}, \ref{D2} that the levels of all the listed structures are equal to $2$.

Let us now prove that $A$ can be represented by an algebra from the corresponding set if $lev(A)=2$. It follows from Lemmas \ref{AnyToSt} and \ref{St} that $G(A)\le 2$. If $G(A)=1$, then $A$ can be represented by a structure from the corresponding set by Proposition \ref{GenType1}. If $G(A)=2$, then there exists a $2$-dimensional $1$-generated subalgebra $C$ of $A$. Let $B$ be the IW contraction of $A$ with respect to $C$. It is clear that $B$ is a trivial singular extension of $C$.  If $lev(B)=1$, then it follows from Proposition \ref{level1} that $C$ is an ideal of $A$ isomorphic to ${\bf A}_3$ as an algebra, and hence the required assertion  follows from Proposition \ref{A3ideal}. If $lev(B)\ge 2$, then $lev(A)=2$ iff $lev(B)=2$ and $A\cong B$, i.e. it remains to consider the case where $A$ is a trivial singular extension of a $2$-dimensional $1$-generated algebra. In this case the required result follows from Lemmas \ref{A3ext}, \ref{nondiag}, \ref{3levext}, \ref{A1}, \ref{A2}, \ref{B2}, \ref{D2}, Corollary \ref{notA3} and the classification of $2$-dimensional $1$-generated algebras.
\end{Proof}

Let us recall that due to \cite{gorb93} the $\infty$-level of an $n$-dimensional algebra $A$ is $lev_{\infty}A=\lim_{m\to\infty}lev_m(A\oplus{\bf k}^{m-n})$. We say that the $n$-dimensional algebra $A$ is {\it stably isomorphic} to the $m$-dimensional algebra $B$ if $A\oplus{\bf k}^{\max(n,m)-n}\cong B\oplus{\bf k}^{\max(n,m)-m}$. It is clear that the $\infty$-level of an algebra is invariant under stable isomorphisms. The next corollary gives the classification of algebras with the $\infty$-level $2$ modulo stable isomorphism, and thus recovers partially the results of \cite{gorb93}, where the anticommutative algebras of the $\infty$-levels $2$ and $3$ were classified. On the other hand, the classification of anticommutative algebras of the $\infty$-level $3$ given in \cite{gorb93} is absolutely wrong, and hence, in fact, we recover all the valid results of this paper. Note also that the classification of algebras with a given $\infty$-level is a much easier problem than the classification of $n$-dimensional algebras with a given level. Some specific methods for this classification are presented in \cite{gorb93}.  Note that it follows from Proposition \ref{level1} that $lev_{\infty}A=1$ iff $A$ is stably isomorphic to an algebra represented by either ${\bf A}_3$ or ${\bf n}_3$.

\begin{corollary} The algebra $A$ has the $\infty$-level $2$ iff it is stably isomorphic to an algebra represented by some structure from the set
$\left\{T_0^{2,\overline{1,0}},T_{0}^{2,2},\eta_2\right\}\cup \{{\bf F}^{\alpha,\beta}\}_{(\alpha,\beta)\in K_{2}^*}.$
\end{corollary}

Finally, at the end of our paper we present corollaries that give the classification of algebras of the level $2$ in some varieties. All of them follow directly from Theorem \ref{ClassLev2}.
In particular, we recover the results of \cite{khud15} for Jordan algebras and correct the results of the same paper for associative algebras.

\begin{corollary} Suppose that ${\rm char}{\bf k}\not=2$. Let $A$ be a commutative $n$-dimensional algebra.
\begin{enumerate}
\item If $n=2$, then $lev(A)=2$ iff $A$ can be represented by some structure from the set $$\left\{{\bf A}_1^{\frac{1}{2}},{\bf B}_2^{\frac{1}{2}}\right\}\cup\{{\bf D}_2^{\alpha,\alpha}\}_{\alpha\in{\bf k}\setminus\left\{\frac{1}{2}\right\}}.$$
\item If $n\ge 3$, then $lev(A)=2$ iff $A$ can be represented by a structure from the set
$$
\left\{{\bf k}^{n-2}\rtimes_{\frac{1}{2}}{\bf A}_1^{\frac{1}{2}},{\bf k}^{n-2}\rtimes_0^t{\bf B}_2^{\frac{1}{2}},{\bf F}^{1,1}\oplus{\bf k}^{n-3}\right\}\cup\{{\bf k}^{n-2}\rtimes_0^t{\bf D}_2^{\alpha,\alpha}\}_{\alpha\in{\bf k}\setminus\left\{\frac{1}{2}\right\}}.
$$
\end{enumerate}
In particular, the set of $n$-dimensional Jordan algebra structures of the level $2$ is formed by the structures ${\bf D}_2^{0,0}$ and ${\bf D}_2^{1,1}$ if $n=2$ and by the structures ${\bf k}^{n-2}\rtimes_0^t{\bf D}_2^{0,0}$, ${\bf k}^{n-2}\rtimes_0^t{\bf D}_2^{1,1}$ and ${\bf F}^{1,1}\oplus{\bf k}^{n-3}$ if $n\ge 3$.
\end{corollary}

Let us recall that the algebra $A$ is called left alternative if $(aa)b=a(ab)$ for all $a,b\in A$. It is clear that an associative algebra is always left alternative.

\begin{corollary} Let $A$ be a left alternative $n$-dimensional algebra.
\begin{enumerate}
\item If $n=2$, then $lev(A)=2$ iff $A$ can be represented either by the structure ${\bf D}_2^{0,0}$ or by the structure ${\bf D}_2^{1,1}$.
\item If $n=3$, then $lev(A)=2$ iff $A$ can be represented by a structure from the set
$$
\left\{{\bf k}\rtimes_0^t{\bf D}_2^{0,0},{\bf k}\rtimes_0^t{\bf D}_2^{1,1},T_1^{2,\overline{1,0}}\right\}\cup\{{\bf F}^{\alpha,\beta}\}_{(\alpha,\beta)\in K_{2}^*}.
$$
\item If $n=4$, then $lev(A)=2$ iff $A$ can be represented by a structure from the set
$$
\left\{{\bf k}^2\rtimes_0^t{\bf D}_2^{0,0},{\bf k}^2\rtimes_0^t{\bf D}_2^{1,1},T_1^{2,\overline{1,0}},T_1^{2,\overline{0,1}}\right\}\cup\{{\bf F}^{\alpha,\beta}\oplus{\bf k}\}_{(\alpha,\beta)\in K_{2}^*}.
$$
\item If $n\ge 5$, then $lev(A)=2$ iff $A$ can be represented by a structure from the set
$$
\left\{{\bf k}^{n-2}\rtimes_0^t{\bf D}_2^{0,0},{\bf k}^{n-2}\rtimes_0^t{\bf D}_2^{1,1},T_1^{2,\overline{1,0}},T_1^{2,\overline{0,1}},T_{0}^{2,2},\eta_2\oplus{\bf k}^{n-5}\right\}\cup\{{\bf F}^{\alpha,\beta}\oplus{\bf k}^{n-3}\}_{(\alpha,\beta)\in K_{2}^*}.
$$
\end{enumerate}
In particular, all the left alternative algebras of the level $2$ are associative.
\end{corollary}

\section{Tables and algebras used in the paper}\label{notat}

\begin{center}
Table 1. {\it Nilpotent algebras of the first $5$ levels in $\mathcal{T}_n$.}\\
\begin{tabular}{|l|l|l|l|}
\hline
{\rm level}&{\rm partition and dimension}&{\rm notation}&{\rm multiplication table}\\
\hline
{\rm 1}&$\begin{array}{l}(2,1,\dots,1),\,\,n\ge 3\end{array}$&$\begin{array}{l}{\bf n}_3\oplus{\bf k}^{n-3}\end{array}$&$\begin{array}{l}e_1e_2=e_3,\,\,e_2e_1=-e_3\end{array}$\\
\hline
\hline
{\rm 2}&$\begin{array}{l}(3),\,\,n=4\vspace{0.1cm}\\(2,2,1,\dots,1),\,\,n\ge 5\end{array}$ &$\begin{array}{l}T_0^3\vspace{0.1cm}\\T_0^{2,2}\end{array}$ & $\begin{array}{l}e_1e_2=e_3,\,\,e_2e_1=-e_3,\,\,e_1e_3=e_4,\,\,e_3e_1=-e_4\vspace{0.1cm}\\e_1e_2=e_3,\,\,e_2e_1=-e_3,\,\,e_1e_4=e_5,\,\,e_4e_1=-e_5\end{array}$\\
\hline
\hline
{\rm 3}&$\begin{array}{l}(3,1,\dots,1),\,\,n\ge 5\vspace{0.1cm}\\(2,2,2,1,\dots,1),\,\,n\ge 7\end{array}$ &$\begin{array}{l}T_0^3\vspace{0.1cm}\\T_0^{2,2,2}\end{array}$ & $\begin{array}{l}e_1e_2=e_3,\,\,e_2e_1=-e_3,\,\,e_1e_3=e_4,\,\,e_3e_1=-e_4\vspace{0.1cm}\\e_1e_{2i}=-e_{2i}e_1=e_{2i+1},\,\,1\le i\le 3\end{array}$\\
\hline
\hline
{\rm 4}&$\begin{array}{l}(4),\,\,n=5\vspace{0.1cm}\\(3,2,1,\dots,1),\,\,n\ge 6\vspace{0.1cm}\\(2,2,2,2,1,\dots,1),\,\,n\ge 9\end{array}$&$\begin{array}{l}T_0^4\vspace{0.1cm}\\T_0^{3,2}\vspace{0.1cm}\\T_0^{2,2,2,2}\end{array}$ & $\begin{array}{l}e_1e_i=e_{i+1},\,\,e_ie_1=-e_{i+1},\,\,2\le i\le 4\vspace{0.1cm}\\e_1e_i=e_{i+1},\,\,e_ie_1=-e_{i+1},\,\,i\in\{2,3,5\}\vspace{0.1cm}\\e_1e_{2i}=e_{2i+1},\,\,e_{2i}e_1=-e_{2i+1},\,\,1\le i\le 4\end{array}$\\
\hline
\hline
{\rm 5}&$\begin{array}{l}(4,1,\dots,1),\,\,n\ge 6\vspace{0.1cm}\\(3,3),\,\,n=7\vspace{0.1cm}\\(3,2,2,1,\dots,1),\,\,n\ge 8\vspace{0.1cm}\\(2,2,2,2,2,1,\dots,1),\,\,n\ge 11\end{array}$ &$\begin{array}{l}T_0^4\vspace{0.1cm}\\T_0^{3,3}\vspace{0.1cm}\\T_0^{3,2,2}\vspace{0.1cm}\\T_0^{2,2,2,2,2}\end{array}$ & $\begin{array}{l}e_1e_i=e_{i+1},\,\,e_ie_1=-e_{i+1},\,\,2\le i\le 4\vspace{0.1cm}\\e_1e_i=-e_ie_1=e_{i+1},\,\,i\in\{2,3,5,6\}\\e_1e_i=-e_ie_1=e_{i+1},\,\,i\in\{2,3,5,7\}\vspace{0.1cm}\\e_1e_{2i}=-e_{2i}e_1=e_{2i+1},\,\,1\le i\le 5\end{array}$\\
\hline
\end{tabular}
\end{center}
\vspace{0.5cm}
\begin{center}
Table 2. {\it Solvable nonnilpotent algebras of the first $5$ levels in $\mathcal{T}_n$.}\\
\begin{tabular}{|l|l|l|l|}
\hline
{\rm level}&{\rm blocks and dimension}&{\rm notation}&{\rm multiplication table}\\
\hline
{\rm 1}&$\begin{array}{l}J(1),\dots,J(1),\,\,n\ge 2\end{array}$&$\begin{array}{l}{\bf p}^-\end{array}$&$\begin{array}{l}e_1e_i=e_i,\,\,e_ie_1=-e_i,\,\,2\le i\le n\end{array}$\\
\hline
\hline
{\rm 2}&$\begin{array}{l}J(\alpha,\beta),\,J(\beta),\dots J(\beta),\,\,n\ge 3,\\(\alpha,\beta)\in\begin{cases}K^*_{2},&\mbox{ if $n=3$,}\\K^*_{1,1},&\mbox{ if $n>3$.}\end{cases}\end{array}$ &
$\begin{array}{l}T_0^{2,\overline{\alpha,\beta}}\end{array}$ & $\begin{array}{l}e_1e_2=\alpha e_2+e_3,\,\,e_2e_1=-\alpha e_2-e_3,\\e_1e_i=\beta e_i,\,\,e_ie_1=-\beta e_i,\,\,3\le i\le n\end{array}$\\
\hline
\end{tabular}
\end{center}
\begin{center}
\begin{tabular}{|l|l|l|l|}
\hline
{\rm level}&{\rm blocks and dimension}&{\rm notation}&{\rm multiplication table}\\
\hline
{\rm 3}&$\begin{array}{l}J(\alpha,\beta,\gamma),\,\,n=4,\,\,(\alpha,\beta,\gamma)\in K^*_3\!\!\!\!\!\!\!\!\!\!\begin{array}{l}\\ \\ \\ \end{array} \vspace{0.2cm}\\
\!\!\!\begin{array}{l}J(\alpha,\beta),\,\,J(\alpha,\beta),\,\,J(\beta),\dots,J(\beta),\,\,n\ge 5,\\(\alpha,\beta)\in\begin{cases}K^*_{2},&\mbox{ if $n=5$,}\\K^*_{1,1},&\mbox{ if $n>5$.}\end{cases}\end{array}\!\!\!\!\!\!\!\!\!\!\begin{array}{l}\\ \\ \\ \\ \end{array} \end{array}$ &$\begin{array}{l}T_0^{3,\overline{\alpha,\beta,\gamma}}\!\!\!\!\!\!\!\!\!\!\begin{array}{l}\\ \\ \\ \end{array}\vspace{0.2cm}\\T_0^{2,2,\overline{\alpha,\beta}}\!\!\!\!\!\!\!\!\!\!\begin{array}{l}\\ \\ \\ \\ \end{array}\end{array}$ & $\begin{array}{l}e_1e_2=\alpha e_2+e_3,\,\,e_2e_1=-\alpha e_2-e_3,\\
e_1e_3=\beta e_3+e_4,\,\,e_3e_1=-\beta e_3-e_4,\\e_1e_4=\gamma e_4,\,\,e_4e_1=-\gamma e_4\vspace{0.2cm}\\
e_1e_2=\alpha e_2+e_3,\,\,e_2e_1=-\alpha e_2-e_3,\\e_1e_3=\beta e_3,\,\,e_3e_1=-\beta e_3,\\e_1e_4=\alpha e_4+e_5,\,\,e_4e_1=-\alpha e_4-e_5,\\e_1e_i=\beta e_i,\,\,e_ie_1=-\beta e_i,\,\,5\le i\le n\end{array}$\\
\hline
\hline
{\rm 4}&$\begin{array}{l}\!\!\!\begin{array}{l}J(\alpha,\beta,\gamma),\,\,J(\gamma),\dots,J(\gamma),\,\,n\ge 5,\\(\alpha,\beta,\gamma)\in K_{2,1}^*\end{array}\!\!\!\!\!\!\!\!\!\!\begin{array}{l}\\ \\ \\ \end{array}\vspace{0.2cm}\\J(\alpha,\beta),\,\,J(\alpha,\beta),\,\,J(\alpha,\beta),\,\,J(\beta),\dots,J(\beta),\\n\ge 7,\,\,(\alpha,\beta)\in\begin{cases}K^*_{2},&\mbox{ if $n=7$,}\\K^*_{1,1},&\mbox{ if $n>7$.}\end{cases} \end{array}$ &
$\begin{array}{l}T_0^{3,\overline{\alpha,\beta,\gamma}}\!\!\!\!\!\!\!\!\!\!\begin{array}{l}\\ \\ \\ \end{array}\vspace{0.2cm}\\T_0^{2,2,2,\overline{\alpha,\beta}}\!\!\!\!\!\!\!\!\!\!\begin{array}{l}\\ \\ \\ \end{array}\end{array}$ &
$\begin{array}{l}e_1e_2=\alpha e_2+e_3,\,\,e_2e_1=-\alpha e_2-e_3,\\e_1e_3=\beta e_3+e_4,\,\,e_3e_1=-\beta e_3-e_4,\\e_1e_i=\gamma e_i,\,\,e_ie_1=-\gamma e_i,\,\,4\le i\le n\vspace{0.2cm}\\
e_1e_i=\alpha e_i+e_{i+1},\,\,e_ie_1=-\alpha e_i-e_{i+1},\\i\in\{2,4,6\},\,\,e_1e_i=\beta e_i,\,\,e_ie_1=-\beta e_i,\\
i\in\{3,5\}\mbox{ \rm and }7\le i\le n\end{array}$\\
\hline
\hline
{\rm 5}&$\begin{array}{l}J(\alpha,\beta,\gamma,\delta),\,\,n=5,\,\,(\alpha,\beta,\gamma,\delta)\in K_4^*\!\!\!\!\!\!\!\!\!\!\begin{array}{l}\\ \\ \\ \\ \end{array}\vspace{0.2cm}\\\!\!\!\begin{array}{l}J(\alpha,\beta,\gamma),\,\,J(\beta,\gamma),\,\,J(\gamma),\dots,J(\gamma),\,\,n\ge 6,\\(\alpha,\beta,\gamma)\in\begin{cases}K^*_{1,2},&\mbox{ if $n=6$,}\\K^*_{1,1,1},&\mbox{ if $n>6$.}\end{cases}\end{array}\!\!\!\!\!\!\!\!\!\!\begin{array}{l}\\ \\ \\ \\ \\ \end{array}\vspace{0.2cm}\\J(\alpha,\beta),\,\,J(\alpha,\beta),\,\,J(\alpha,\beta),\,\,J(\alpha,\beta),\\
J(\beta),\dots,J(\beta),\,\,n\ge 9,\\(\alpha,\beta)\in\begin{cases}K^*_{2},&\mbox{ if $n=9$,}\\K^*_{1,1},&\mbox{ if $n>9$.}\end{cases}\end{array}$ & $\begin{array}{l}T_0^{4,\overline{\alpha,\beta,\gamma,\delta}}\!\!\!\!\!\!\!\!\!\!\begin{array}{l}\\ \\ \\ \\ \end{array}\vspace{0.2cm}\\T_0^{3,2,\overline{\alpha,\beta,\gamma}}\!\!\!\!\!\!\!\!\!\!\begin{array}{l}\\ \\ \\ \\ \\ \end{array}\vspace{0.2cm}\\T_0^{2,2,2,2,\overline{\alpha,\beta}}\!\!\!\!\!\!\!\!\!\!\begin{array}{l}\\ \\ \\ \\ \end{array}\end{array}$ &
$\begin{array}{l}e_1e_2=\alpha e_2+e_3,\,\,e_2e_1=-\alpha e_2-e_3,\\e_1e_3=\beta e_3+e_4,\,\,e_3e_1=-\beta e_3-e_4,\\e_1e_4=\gamma e_4+e_5,\,\,e_4e_1=-\gamma e_4-e_5,\\e_1e_i=\delta e_i,\,\,e_ie_1=-\delta e_i,\,\,5\le i\le n\vspace{0.2cm}\\
e_1e_2=\alpha e_2+e_3,\,\,e_2e_1=-\alpha e_2-e_3,\\e_1e_3=\beta e_3+e_4,\,\,e_3e_1=-\beta e_3-e_4,\\e_1e_4=\gamma e_4,\,\,e_4e_1=-\gamma e_4,\\e_1e_5=\beta e_5+e_6,\,\,e_5e_1=-\beta e_5-e_6,\\e_1e_i=\gamma e_i,\,\,e_ie_1=-\gamma e_i,\,\,6\le i\le n\vspace{0.2cm}\\
\!\!\!\begin{array}{l}e_1e_i=\alpha e_i+e_{i+1},\,\,e_ie_1=-\alpha e_i-e_{i+1},\\i\in\{2,4,6,8\},\,\,e_1e_i=\beta e_i,\,\,e_ie_1=-\beta e_i,\\
i\in\{3,5,7\}\mbox{ \rm and }9\le i\le n\end{array}\!\!\!\!\!\!\!\!\!\!\begin{array}{l}\\ \\ \\ \\ \end{array}\end{array}$\\
\hline
\end{tabular}
\end{center}
\vspace{0.5cm}
\begin{center}
Table 3. {\it Nonsolvable algebras of the first $5$ levels in $\mathcal{T}_n$.}\\
\begin{tabular}{|l|l|l|l|}
\hline
{\rm level}&{\rm blocks and dimension}&{\rm notation}&{\rm multiplication table}\\
\hline
{\rm 1}&$\begin{array}{l}J(\alpha),\dots,J(\alpha),\,\,n\ge 2,\,\,\alpha\in{\bf k}\begin{array}{l}\\ \\\end{array}\end{array}$&$\begin{array}{l}\nu^{\alpha}\begin{array}{l}\\ \\\end{array}\end{array}$&$\begin{array}{l}e_1e_1=e_1,\\e_1e_i=\alpha e_i,\,\,e_ie_1=(1-\alpha)e_i,\,\,2\le i\le n\end{array}$\\
\hline
\hline
{\rm 2}&$\begin{array}{l}J(\alpha,\beta),\,J(\beta),\dots J(\beta),\,\,n\ge 3,\\(\alpha,\beta)\in\begin{cases}K_{2},&\mbox{ if $n=3$,}\\K_{1,1},&\mbox{ if $n>3$.}\end{cases}\end{array}$ &
$\begin{array}{l}T_1^{2,\overline{\alpha,\beta}}\end{array}$ & $\begin{array}{l}e_1e_1=e_1,\,\,e_1e_2=\alpha e_2+e_3,\\e_2e_1=(1-\alpha) e_2-e_3,\\e_1e_i=\beta e_i,\,\,e_ie_1=(1-\beta) e_i,\,\,3\le i\le n\end{array}$\\
\hline
\end{tabular}
\end{center}
\begin{center}
\begin{tabular}{|l|l|l|l|}
\hline
{\rm level}&{\rm blocks and dimension}&{\rm notation}&{\rm multiplication table}\\
\hline
{\rm 3}&$\begin{array}{l}J(\alpha,\beta,\gamma),\,\,n=4,\,\,(\alpha,\beta,\gamma)\in K_3\!\!\!\!\!\!\!\!\!\!\begin{array}{l}\\ \\ \\ \\ \end{array} \vspace{0.2cm}\\
\!\!\!\begin{array}{l}J(\alpha,\beta),\,\,J(\alpha,\beta),\,\,J(\beta),\dots,J(\beta),\,\,n\ge 5,\\(\alpha,\beta)\in\begin{cases}K_{2},&\mbox{ if $n=5$,}\\K_{1,1},&\mbox{ if $n>5$.}\end{cases}\end{array}\!\!\!\!\!\!\!\!\!\!\begin{array}{l}\\ \\ \\ \\ \\ \end{array} \end{array}$ &$\begin{array}{l}T_1^{3,\overline{\alpha,\beta,\gamma}}\!\!\!\!\!\!\!\!\!\!\begin{array}{l}\\ \\ \\ \\ \end{array}\vspace{0.2cm}\\T_1^{2,2,\overline{\alpha,\beta}}\!\!\!\!\!\!\!\!\!\!\begin{array}{l}\\ \\ \\ \\ \\ \end{array}\end{array}$ & $\begin{array}{l}e_1e_1=e_1,\,\,e_1e_2=\alpha e_2+e_3,\\e_2e_1=(1-\alpha) e_2-e_3,\\
e_1e_3=\beta e_3+e_4,\,\,e_3e_1=(1-\beta) e_3-e_4,\\e_1e_4=\gamma e_4,\,\,e_4e_1=(1-\gamma) e_4\vspace{0.2cm}\\
e_1e_1=e_1,\,\,e_1e_2=\alpha e_2+e_3,\\e_2e_1=(1-\alpha) e_2-e_3,\\e_1e_3=\beta e_3,\,\,e_3e_1=(1-\beta) e_3,\\e_1e_4=\alpha e_4+e_5,\,\,e_4e_1=(1-\alpha) e_4-e_5,\\e_1e_i=\beta e_i,\,\,e_ie_1=(1-\beta) e_i,\,\,5\le i\le n\end{array}$\\
\hline
\hline
{\rm 4}&$\begin{array}{l}\!\!\!\begin{array}{l}J(\alpha,\beta,\gamma),\,\,J(\gamma),\dots,J(\gamma),\,\,n\ge 5,\\(\alpha,\beta,\gamma)\in K_{2,1}\end{array}\!\!\!\!\!\!\!\!\!\!\begin{array}{l}\\ \\ \\ \\ \end{array}\vspace{0.2cm}\\J(\alpha,\beta),\,\,J(\alpha,\beta),\,\,J(\alpha,\beta),\,\,J(\beta),\dots,J(\beta),\\n\ge 7,\,\,(\alpha,\beta)\in\begin{cases}K_{2},&\mbox{ if $n=7$,}\\K_{1,1},&\mbox{ if $n>7$.}\end{cases} \end{array}$ &
$\begin{array}{l}T_1^{3,\overline{\alpha,\beta,\gamma}}\!\!\!\!\!\!\!\!\!\!\begin{array}{l}\\ \\ \\ \\ \end{array}\vspace{0.2cm}\\T_1^{2,2,2,\overline{\alpha,\beta}}\!\!\!\!\!\!\!\!\!\!\begin{array}{l}\\ \\ \\ \\ \end{array}\end{array}$ &
$\begin{array}{l}e_1e_1=e_1,\,\,e_1e_2=\alpha e_2+e_3,\\e_2e_1=(1-\alpha) e_2-e_3,\\e_1e_3=\beta e_3+e_4,\,\,e_3e_1=(1-\beta) e_3-e_4,\\e_1e_i=\gamma e_i,\,\,e_ie_1=(1-\gamma) e_i,\,\,4\le i\le n\vspace{0.2cm}\\
e_1e_1=e_1,\,\,e_1e_i=\alpha e_i+e_{i+1},\\e_ie_1=(1-\alpha) e_i-e_{i+1},i\in\{2,4,6\},\\e_1e_i=\beta e_i,\,\,e_ie_1=(1-\beta) e_i,\\
i\in\{3,5\}\mbox{ \rm and }7\le i\le n\end{array}$\\
\hline
\hline
{\rm 5}&$\begin{array}{l}J(\alpha,\beta,\gamma,\delta),\,\,n=5,\,\,(\alpha,\beta,\gamma,\delta)\in K_4\!\!\!\!\!\!\!\!\!\!\begin{array}{l}\\ \\ \\ \\ \\ \end{array}\vspace{0.2cm}\\\!\!\!\begin{array}{l}J(\alpha,\beta,\gamma),\,\,J(\beta,\gamma),\,\,J(\gamma),\dots,J(\gamma),\,\,n\ge 6,\\(\alpha,\beta,\gamma)\in\begin{cases}K_{1,2},&\mbox{ if $n=6$,}\\K_{1,1,1},&\mbox{ if $n>6$.}\end{cases}\end{array}\!\!\!\!\!\!\!\!\!\!\begin{array}{l}\\ \\ \\ \\ \\ \\ \end{array}\vspace{0.2cm}\\J(\alpha,\beta),\,\,J(\alpha,\beta),\,\,J(\alpha,\beta),\,\,J(\alpha,\beta),\\
J(\beta),\dots,J(\beta),\,\,n\ge 9,\\(\alpha,\beta)\in\begin{cases}K_{2},&\mbox{ if $n=9$,}\\K_{1,1},&\mbox{ if $n>9$.}\end{cases}\vspace{0.15cm}\end{array}$ & $\begin{array}{l}T_1^{4,\overline{\alpha,\beta,\gamma,\delta}}\!\!\!\!\!\!\!\!\!\!\begin{array}{l}\\ \\ \\ \\ \\ \end{array}\vspace{0.2cm}\\T_1^{3,2,\overline{\alpha,\beta,\gamma}}\!\!\!\!\!\!\!\!\!\!\begin{array}{l}\\ \\ \\ \\ \\ \\ \end{array}\vspace{0.2cm}\\T_1^{2,2,2,2,\overline{\alpha,\beta}}\!\!\!\!\!\!\!\!\!\!\begin{array}{l}\\ \\ \\ \\ \\ \end{array}\end{array}$ &
$\begin{array}{l}e_1e_1=e_1,\,\,e_1e_2=\alpha e_2+e_3,\\e_2e_1=(1-\alpha) e_2-e_3,\\e_1e_3=\beta e_3+e_4,\,\,e_3e_1=(1-\beta) e_3-e_4,\\e_1e_4=\gamma e_4+e_5,\,\,e_4e_1=(1-\gamma) e_4-e_5,\\e_1e_i=\delta e_i,\,\,e_ie_1=(1-\delta) e_i,\,\,5\le i\le n\vspace{0.2cm}\\
e_1e_1=e_1,\,\,e_1e_2=\alpha e_2+e_3,\\e_2e_1=(1-\alpha) e_2-e_3,\\e_1e_3=\beta e_3+e_4,\,\,e_3e_1=(1-\beta) e_3-e_4,\\e_1e_4=\gamma e_4,\,\,e_4e_1=(1-\gamma) e_4,\\e_1e_5=\beta e_5+e_6,\,\,e_5e_1=(1-\beta0 e_5-e_6,\\e_1e_i=\gamma e_i,\,\,e_ie_1=(1-\gamma) e_i,\,\,6\le i\le n\vspace{0.2cm}\\
\!\!\!\begin{array}{l}e_1e_1=e_1,\,\,e_1e_i=\alpha e_i+e_{i+1},\\e_ie_1=(1-\alpha) e_i-e_{i+1},\,\,i\in\{2,4,6,8\},\\e_1e_i=\beta e_i,\,\,e_ie_1=(1-\beta) e_i,\\
i\in\{3,5,7\}\mbox{ \rm and }9\le i\le n\end{array}\!\!\!\!\!\!\!\!\!\!\begin{array}{l}\\ \\ \\ \\ \\ \end{array}\end{array}$\\
\hline
\end{tabular}
\end{center}
\vspace{0.5cm}
\newpage
\begin{center}
Table 4. {\it $2$-dimensional algebras.}\\
\begin{tabular}{|l|l|}
\hline
${\bf A}_1^{\alpha}$, $\alpha\in{\bf k}$ & $e_1e_1=e_1+e_2,$ $e_1e_2=\alpha e_2,$ $e_2e_1= (1-\alpha) e_2$ \\

\hline
${\bf A}_2$ & $e_1e_1=e_2,$ $e_1e_2=e_2,$ $e_2e_1= -e_2$ \\

\hline
${\bf A}_3$ & $e_1e_1=e_2$\\

\hline
${\bf A}_4^{\alpha}$, $\alpha\in{\bf k}$ & $e_1e_1= \alpha e_1+  e_2,$ $e_1e_2= e_1+ \alpha e_2,$ $ e_2e_1= - e_1$\\

\hline
${\bf B}_1^{\alpha}$, $\alpha\in{\bf k}$ &  $e_1e_2=(1-\alpha)e_1 +  e_2,$ $ e_2e_1=\alpha  e_1 - e_2$\\

\hline
${\bf B}_2^{\alpha}$, $\alpha\in{\bf k}$ &  $e_1e_2=\alpha e_2,$ $ e_2e_1=(1-\alpha)  e_2$\\

\hline
${\bf C}^{\alpha,\beta}$, $\alpha,\beta \in{\bf k}$ & $e_1e_1=e_2,$ $ e_1e_2= (1-\alpha)e_1+\beta e_2,$ $e_2e_1=\alpha e_1-\beta e_2,$ $ e_2e_2=e_2$ \\

\hline
${\bf D}_1^{\alpha,\beta}$, $\alpha,\beta\in {\bf k}$ & $e_1e_1=e_1,$ $ e_1e_2=(1-\alpha)e_1+ \beta  e_2,$ $e_2e_1=\alpha  e_1 - \beta   e_2$\\

\hline
${\bf D}_{2}^{\alpha,\beta}$, $\alpha,\beta\in{\bf k}$, $\alpha+\beta\not=1$ & $e_1e_1=e_1,$ $ e_1e_2=\alpha e_2,$ $e_2e_1=\beta e_2$\\

\hline
${\bf D}_{3}^{\alpha,\beta}$, $\alpha,\beta\in{\bf k}$, $\alpha+\beta\not=1$ & $e_1e_1=e_1,$ $ e_1e_2=e_1+ \alpha  e_2,$ $e_2e_1=- e_1 + \beta  e_2$\\

\hline
${\bf E}_1^{\alpha,\beta,\gamma,\delta}$, $\alpha,\beta,\gamma,\delta\in {\bf k}$, $\beta+\delta\not=1$ & $e_1e_1=e_1,$ $ e_1e_2=\alpha e_1+ \beta  e_2,$ $e_2e_1=\gamma  e_1 + \delta  e_2,$ $ e_2e_2=e_2$\\

\hline
${\bf E}_4$ & $e_1e_1=e_1,$ $ e_1e_2=e_1+  e_2,$ $ e_2e_2=e_2$\\

\hline
\end{tabular}
\end{center}
\vspace{0.5cm}
\begin{center}
Table 5. {\it Trivial singular extensions of $2$-dimensional algebras.}\\
\begin{tabular}{|l|l|}
\hline
$\begin{array}{l}{\bf k}^{n-2}\rtimes_{\epsilon} {\bf A}_1^{\alpha},\,\alpha,\epsilon\in{\bf k}\end{array}$ & $\begin{array}{l}e_1e_1=e_1+e_2,\,e_1e_2=\alpha e_2,\,e_2e_1=(1-\alpha)e_2,\,e_1e_i=\epsilon e_i,\,e_ie_1=(1-\epsilon)e_i\,\,3\le i\le n\end{array}$ \\

\hline
$\begin{array}{l}{\bf k}^{n-2}\rtimes_{\epsilon} {\bf A}_2,\,\epsilon\in{\bf k}\end{array}$ & $\begin{array}{l}e_1e_1=e_2,\,e_1e_2=e_2,\,e_2e_1=-e_2,\,e_1e_i=\epsilon e_i,\,e_ie_1=-\epsilon e_i\,\,3\le i\le n\end{array}$\\

\hline

$\begin{array}{l}{\bf k}^{n-2}\rtimes{\bf A}_3\end{array}$ & $\begin{array}{l}e_1e_1=e_2,\,e_1e_i=e_i,\,e_ie_1=-e_i,\,\,3\le i\le n\end{array}$\\

\hline

$\begin{array}{l}{\bf k}^2\rtimes_{\alpha,\beta}{\bf A}_3\,\,(n=4),\\(\alpha,\beta)\in K_{1,1}^*\end{array}$ & $\begin{array}{l}e_1e_1=e_2,\,e_1e_3=\alpha e_4,\,e_4e_1=\beta e_4\end{array}$\\

\hline

$\begin{array}{l}{\bf k}^{n-2}\rtimes_{\epsilon} {\bf B}_2^{\alpha},\,\alpha,\epsilon\in{\bf k}\end{array}$ & $\begin{array}{l}e_1e_2=\alpha e_2,\,e_2e_1=(1-\alpha)e_2,\,e_1e_i=\epsilon e_i,\,e_ie_1=-\epsilon e_i\,\,3\le i\le n\end{array}$\\

\hline

$\begin{array}{l}{\bf k}^{n-2}\rtimes_{\epsilon}^t {\bf B}_2^{\alpha},\\\alpha\in{\bf k},\,\epsilon\in\{0,1\}\end{array}$ & $\begin{array}{l}e_1e_i=\alpha e_i,\,e_ie_1=(1-\alpha)e_i\,\,2\le i\le n,\,e_2e_i=\epsilon e_i,\,e_ie_2=-\epsilon e_i\,\,3\le i\le n\end{array}$\\

\hline

$\begin{array}{l}{\bf k}^{n-2}\rtimes_{\epsilon} {\bf D}_2^{\alpha,\beta},\\\alpha,\beta,\epsilon\in{\bf k},\,\alpha+\beta\not=1\end{array}$ & $\begin{array}{l}e_1e_1=e_1,\,e_1e_2=\alpha e_2,\,e_2e_1=\beta e_2,\,e_1e_i=\epsilon e_i,\,e_ie_1=(1-\epsilon) e_i\,\,3\le i\le n\end{array}$\\

\hline

$\begin{array}{l}{\bf k}^{n-2}\rtimes_{\epsilon}^t {\bf D}_2^{\alpha,\beta},\,\alpha,\beta\in{\bf k},\\\epsilon\in\{0,1\},\,\alpha+\beta\not=1\end{array}$ & $\begin{array}{l}e_1e_1=e_1,\,e_1e_i=\alpha e_i,\,e_ie_1=\beta e_i\,\,2\le i\le n,\,e_2e_i=\epsilon e_i,\,e_ie_2=-\epsilon e_i\,\,3\le i\le n\end{array}$\\

\hline

$\begin{array}{l}{\bf k}^{n-2}\rtimes {\bf E}_4\end{array}$ & $\begin{array}{l}e_1e_1=e_1,\,e_1e_2=e_1+e_2,\,e_2e_2=e_2,\,e_1e_i=e_ie_2=e_i\,\,3\le i\le n\end{array}$\\

\hline
\end{tabular}
\end{center}
\vspace{0.5cm}
\begin{center}
Table 6. {\it Other algebra structures used in the paper.}\\
\begin{tabular}{|ll|l|}
\hline
$\begin{array}{l}{\bf G}\,\,\end{array}$ & $n=3$ & $\begin{array}{l}e_1e_1=e_2,\,e_2e_2= e_3\end{array}$\\

\hline
$\begin{array}{l}{\bf G}^{\alpha,\beta}\,\,,\,(\alpha,\beta)\in K_{1,1}^*\end{array}$ & $n=3$ & $\begin{array}{l}e_1e_1=e_2,\,e_1e_2=\alpha e_3,\,e_2e_1=\beta e_3\end{array}$\\

\hline

$\begin{array}{l}\eta_m\,\,\end{array}$ & $n=2m+1$ & $\begin{array}{l}e_{2i-1}e_{2i}=e_{2m+1},\,e_{2i}e_{2i-1}=-e_{2m+1},\,\,1\le i\le m\end{array}$\\

\hline
$\begin{array}{l}{\bf F}^{\alpha,\beta}\in\mathcal{A}_3\,\,\,(\alpha,\beta)\in K_{2}^*\end{array}$ & $n=3$ & $\begin{array}{l}e_1e_1=e_3,\,e_1e_2=\alpha e_3,\,e_2e_1=\beta e_3\end{array}$\\

\hline
$\begin{array}{l}{\bf A}_3\rtimes{\bf p}^-\,\,\end{array}$  & $n\ge 3$ & $\begin{array}{l}e_{n-1}e_{n-1}=e_n,\,e_1e_i=e_i,\,e_ie_1=-e_i,\,\,2\le i\le n\end{array}$ \\

\hline
$\begin{array}{l}{\bf A}_3\rtimes \nu^{\alpha}\,\,\,\alpha\in{\bf k}\end{array}$ & $n\ge 3$ & $\begin{array}{l}e_1e_1=e_1,\,e_{n-1}e_{n-1}=e_n,\,e_1e_i=\alpha e_i,\,e_ie_1=(1-\alpha)e_i,\,\,2\le i\le n\end{array}$\\

\hline

$\begin{array}{l}({\bf A}_3\oplus{\bf k}^{n-4})\rtimes {\bf E}_4\,\,\end{array}$ & $n\ge 4$ & $\begin{array}{l}e_1e_1=e_1,\,e_1e_2=e_1+e_2,\,e_2e_2=e_2,\,e_{n-1}e_{n-1}=e_n,\\e_1e_i=e_ie_2=e_i\,\,3\le i\le n\end{array}$\\

\hline
\end{tabular}
\end{center}

\end{document}